\documentclass[reqno]{amsart}

\usepackage{amsmath}
\usepackage{amsthm}
\usepackage{amssymb}
\usepackage{url}

\usepackage{tikz-cd}
\usepackage[english]{babel}
\usepackage{amsmath,amsthm,amssymb,amscd,amsfonts,mathrsfs,tikz,bm}
\usepackage{xcolor}
\usepackage[latin1]{inputenc}
\usepackage{soul}
\usepackage[all]{xy}
\usetikzlibrary{matrix,arrows,positioning,calc}
\usepackage[T1]{fontenc}
\usepackage{textcomp}
\usepackage{palatino,helvet}
\usepackage{titlesec}
\usepackage[hang,flushmargin]{footmisc}
\usepackage{enumerate} 

\titlelabel{\thetitle.\quad}

\titleformat{\subsection}[runin]{\normalfont\bfseries}{\thesubsection.}{3pt}{}
\titleformat{\subsubsection}[runin]{\normalfont\bfseries}{\thesubsubsection.}{3pt}{}

\numberwithin{equation}{subsection}\theoremstyle{plain}

\theoremstyle{definition}

\theoremstyle{remark}

\def\Z{\mathbb Z}
\def\mcF{\mathcal F}
\def\mcE{\mathcal{E}}
\def\mcC{\mathcal{C}}

\def\mcG{\mathcal G}

\def\Mod{\mbox{-Mod}}

\def\bfC+{\mathbf{C}_+}

\def\C{{\mathcal{C}}}

\def\bfD{\mathbf{D}}
\def\A{\mathcal{A}}
\def\bfC{\mathbf{C}}

\def\A{\mathcal{A}}

\newcommand{\Ext}{\operatorname{Ext}}
\newcommand{\End}{\mbox{End}}

\newcommand{\Add}{\mbox{Add}}

\newcommand{\Hom}{\operatorname{Hom}}

\newcommand{\id}{\operatorname{id}}

\newcommand{\Proj}{\mbox{Proj}}

\newcommand{\coker}{\mbox{coker}}

\newcommand{\Ab}{\operatorname{Ab}}

\newcommand{\ev}{\operatorname{ev}}

\begin{document}

\title[]{Transfer of homological objects in exact categories  via adjoint triples. Applications to functor categories}

\author[S. Estrada]{Sergio Estrada}
\address{
Universidad de Murcia \\
Facultad de Matem\'aticas\\
Spain
}
\email[S. Estrada]{sestrada@um.es}

\author[M. Cort\'es-Izurdiaga]{Manuel Cort\'es-Izurdiaga}
\address{
Universidad de Malaga \\
Facultad de Matem\'atica Aplicada\\
Spain
}
\email[M. Cort\'es-Izurdiaga]{ mizurdiaga@uma.es}

\author[S. Odaba\c si]{S\.{I}nem Odaba\c si}
\address{
Universidad de Murcia \\
Facultad de Matem\'aticas \\
Spain
}
\email[S. Odaba\c si]{sinem.odabasi@um.es}

\date{}

\thanks{2020 Mathematics Subject Classification. 18G35,  18G25, 18E99.}

\maketitle

\begin{abstract}

For a given  family  $\{({ \rm q}_i, { \rm t}_i, { \rm p_i})\}_{i \in I}$  of adjoint triples between  exact categories $\mathcal{C}$ or  $\mathcal{D}$, we show that 
  any cotorsion pair   in $\mathcal{C}$ and $\mathcal{D}$ yield two canonical cotorsion pairs  providing a concrete description of objects   without using any injectives/projectives object hypothesis. We firstly apply this result for the evaluation functor on the functor category $\Add(\A, R\Mod)$  equipped with an exact structure $\mathcal{E}$.  Under mild conditions on $\A$, we introduce the stalk functor at any object of $\A$, and subsequently, we investigate cotorsion pairs induced by stalk functors. Finally, we use them to present an intrinsic characterization of projective/injective objects in $(\Add(\A, R\Mod); \mathcal{E})$.

\end{abstract}


\section*{Introduction}
The lifting of particular homological objects to functor categories has  been of interest for a long time as  it provides a kind of control for globally defined objects through their local information. The emblematic examples come up in  the category of chain complexes. The author in \cite{Gil04} shows that  a complete hereditary cotorsion pair $(\mathcal{F}, \mathcal{G})$ in $R \Mod$, the category of left $R$-modules, induces a pair of hereditary complete cotorsion pairs $(\mbox{dg}\ \tilde{\mathcal{F}}, \tilde{\mathcal{G}} )$ and $(\tilde{\mathcal{F}}, \mbox{dg}\  \tilde{\mathcal{G}} )$ in $\bfC(R)$, the category of chain complexes of left $R$-modules, and provides an intrinsic characterization of chain complexes in $\tilde{\mathcal{F}}$ and $\tilde{\mathcal{G}}$.  Combining it with the well-known result in \cite{Hov02}, those induced cotorsion pairs yield a model structure on $\bfC(R)$ whose homotopy category is just the derived category $\bfD(R)$. There are several advantages of this lifting mechanism. Firstly,   it generates several non-trivial model structures on $\bfC(R)$, as well as recovering the well-known injective and projective ones. Secondly,  the classes $\tilde{\mathcal{F}}$ and $\tilde{\mathcal{G}}$ are mostly the classes of  categorically defined homological objects in $\bfC(R)$ such as projective, injective, flat and absolutely pure chain complexes.   Those results have been  generalized to a category of chain complexes $\bfC(\mathcal{C})$ over a sufficiently nice exact category $\mathcal{C}$; see \cite[Section 7]{Sto13}. In the recent work \cite{HJ22}, the authors focus on the problem of how to obtain projective and injective model structures on $k\mbox{-Lin}( \mathcal{A}, R \Mod)$, the category of $R$-module valued $k$-linear functors on a more general $k$-linear category $\A$, where $k$ is a commutative ring and $R$ is a $k$-algebra. 

Several authors have applied a similar lifting process to modules over a formal triangular matrix ring $L:=\begin{bmatrix}
T & 0\\
N & S
\end{bmatrix}$. For instance, the author in \cite{Mao20} lifts two  cotorsion pairs $(\mcF_T, \mcG_T)$ and $(\mcF_S, \mcG_S)$ in $T\Mod $ and $S \Mod$, respectively, to several cotorsion pairs in  $L\Mod$ providing an explicit description of their modules.  As a  particular case, the trivial cotorsion pairs $(\mcF_T, \mcG_T)$ and $(\mcF_S, \mcG_S)$   (i.e., projective/injective modules) lead to the trivial cotorsion pairs in $L \Mod$. In this case,  the characterization of projective/injective left $L$-modules is the same as in \cite{HV99} and \cite{HV00}.  The subtle difference with  the case of chain complexes  is that the classes used to lift come from the module categories over two rings $T$ and $S$. However, the standard tool used in both cases is the adjoint functors. So
this paper aims to provide a systematic study of how to lift  cotorsion pairs through general adjoint functors between exact categories with applications in functor categories so that it covers the categories of chain complexes and modules over  a formal triangular matrix ring.  Specifically, we give a very intrinsic characterization of projective and injective additive functors under mild conditions.

\sloppy The paper starts by proposing a general question: For a given integer $n \geq 1$, and a pair  $({ \rm q}, { \rm t})$ of adjoint functors between  exact categories $\mathcal{C}$ and $\mathcal{D}$, under which conditions does there exist a canonical natural isomorphism ${\Ext^n({ \rm q}(-),-) \cong \Ext^n(-, { \rm t}(-))}$? Of course, this isomorphism is immediate when both functors are exact between exact categories. However, this is not the case in most situations of interest (for instance, see Example~ \ref{ex:chain_stalk}). In Section 2, we address this question in full generality, fixing one of the objects in components, that is, for given objects $C \in \bfC$ and $D \in \mathcal{D}$, to guarantee  natural isomorphisms 
\begin{equation}\tag{$\ast$}\label{prob:ext}
{\Ext^n({ \rm q}(C),-) \cong \Ext^n(C, { \rm t}(-))}\quad  \textrm{ and  }\quad  {\Ext^n({ \rm q}(-),D) \cong \Ext^n(-, { \rm t}(D))}.
\end{equation}

 The case $n=1$ is fully answered in Proposition~\ref{prop:adjoints}. For it, we introduce two important notions: ${ \rm q}$-flat and ${\rm t}$-coflat objects (see Definition~\ref{def:flat-coflat}). These notions are crucial since they help remove the hypothesis of having enough projective/injective objects required in most situations with  similar constructions such as \cite{HJ19} and \cite{Mao20} (see also Remark~\ref{remark:clases-HJ}). We also show that  the notion of acyclicity for   chain complexes over any exact category can be expressed in terms of (co)-flatness to specific functors (see Example~\ref{ex:chain_stalk}). 
 
When working with a general $n \geq 1$, lengthy  technical conditions  will be required  if one aims to avoid projective/injective objects. Instead,  we firstly rethink the  problem \eqref{prob:ext} not for  fixed objects, but for  fixed classes $\mcF \subseteq \mathcal{C}$ and $\mathcal{G} \subseteq \mathcal{D}$. So under conditions given in Lemma~\ref{prop:heredit1}, we prove that there exist the following  natural isomorphisms
 $${\Ext^n({ \rm q}( \mathcal{F}),-) \cong \Ext^n(\mathcal{F}, { \rm t}(-))}\quad  \textrm{ and  }\quad  {\Ext^n({ \rm q}(-),\mathcal{G}) \cong \Ext^n(-, { \rm t}(\mcG))}.$$
This result is meritorious on its own. The classic example representing this phenomenon is the adjoint triple induced from sphere chain complexes and acyclic chain complexes (see Example~\ref{ex:chain_stalk}).  Lemma~\ref{prop:heredit1}  also provides us the tool of the relatively deriving for functors ${ \rm q}$ and ${\rm t}$ to answer the question  \eqref{prob:ext} for any $n\geq 1$ (see Proposition~\ref{lemma:derived_zero}). 

In Section 3, we apply results from Section~2 for the lifting problem of cotorsion pairs through adjoint functors. Corollary \ref{setup:3}, a special case of Theorem~\ref{Prop:degrewise-cot.pair}, is the main form of the situations of interest, and is used frequently in the rest of the paper. To explain, if there is  a family  $\{({ \rm q}_i, { \rm t}_i, { \rm p_i})\}_{i \in I}$  of adjoint triples between  exact categories
	$$
\xymatrix{\mathcal{C} \ar@/^2em/[rrr]^{{ \rm q}_i} \ar@/_2em/[rrr]_{ { \rm p_i}} &&& \mathcal{D} \ar[lll]^{{ \rm t}_i} }, 
$$
then   any cotorsion pairs $(\mcF, \mcG)$ and $(\mcF', \mcG')$  in $\mathcal{C}$ and $\mathcal{D}$, respectively, yield two canonical cotorsion pairs 
$$({}^\perp( { \rm q}_*(\mcF)^\perp), { \rm q}_*(\mcF)^\perp) \quad \textrm{ and  } \quad ({}^\perp{ \rm p}_*(\mcG),({}^\perp{ \rm p}_*(\mcG))^\perp )$$
$$({}^\perp( { \rm t}_*(\mcF')^\perp), { \rm t}_*(\mcF')^\perp) \quad \textrm{ and  } \quad ({}^\perp{ \rm t}_*(\mcG'),({}^\perp{ \rm t}_*(\mcG'))^\perp )$$
in  $\mathcal{D}$ and  $\mathcal{C}$, respectively. In Corollary~\ref{setup:3}, we  provide a concrete description of objects in ${ \rm q}_*(\mcF)^\perp$, ${}^\perp{ \rm p}_*(\mcG)$, ${ \rm t}_*(\mcF')^\perp$ and  ${}^\perp{ \rm t}_*(\mcG')$  without using any injectives/projectives objects. Applications of this result are presented in the next sections.

Starting from Section~4, we focus on  $\Add(\A, R\Mod)$, the category of $R$-module-valued additive functors on an $R$-$R$-bimodule enriched small category $\A$ (see Hypothesis \ref{hypot_A}).  Our approach is based on consideration of module categories attached locally to $\A$, that is,  the module category over the endomorphism ring $R_A:=\End_{\A}(A)$  of any object $A$ in $\A$.   We assume that the product category $\prod_{A \in \A} R_A \Mod$ is equipped with a  fixed  exact sructure $\prod_{A \in \A} \mathcal{E}_A$, and we consider  a cotorsion pair $(\prod_{A \in \A} \mcF_A, \prod_{A \in \A} \mcG_A)$ in $\prod_{A \in \A} R_A \Mod$.

In Section~4, we also  study  the adjoint triple induced from the evaluation functor $\ev_A: \Add(\A,R\Mod) \rightarrow R \Mod$, $\ev_A(X):=X(A)$. As $\A$ is an $R$-$R$-bimodule enriched category, $\ev_A$ factors through   $R_A \Mod$. So the fixed  exact sructure $\prod_{A \in \A} \mathcal{E}_A$ on the product category $\prod_{A \in \A} R_A \Mod$   induces an  exact structure $\mathcal{E}$  on  $\Add(\A, R\Mod)$ defined as follows: $\mathbb{E} \in \mathcal{E}$ if and only in $\ev_A (\mathbb{E}) \in \mathcal{E}_A$ for every $A \in \A$ (see Hypothesis~\ref{degreewise_exact}). So it makes the functor ${\ev:=(\ev_A)_{A \in \A}:\Add(\A, R\Mod) \longrightarrow  \prod_{A \in \A} R_A \Mod}$ an exact functor. Applying Corollary~\ref{setup:3} for the adjoint triple $({\rm q}, \ev, { \rm p})$, we have :

\vspace{5mm}

\textbf{Proposition~\ref{prop:cotorsion1}:} If for every $A \in \A$, any object in $\mcF_A$ and $\mcG_A$ is ${ \rm q}_A$-flat and ${ \rm p}_A$-coflat, respectively, then  
\begin{align*}
{ \rm q}(  \prod_{A \in \A} \mcF_A)^\perp&=\{ X \in \Add(\mathcal{A}, R\Mod) \mid\ X(A) \in \mcG_A \textrm{ for every }A\in \A\},\\
^\perp{ \rm p}(  \prod_{A \in \A} \mcG_A)&=\{ X \in \Add(\mathcal{A}, R\Mod) \mid\ X(A) \in \mcF_A \textrm{ for every }A\in \A \}.
\end{align*}
	
	\vspace{4mm}
In Section~5, we impose further conditions (see Hypothesis~\ref{hyp2}) on the category $\A$  in order to define ${ \rm s}_A: R_A \Mod \longrightarrow \Add(\A,R\Mod)$, the `stalk' functor at an object $A$ in $\A$ (see \eqref{def:stalk} and Proposition~\ref{prop:stalk1}).  Our stalk functor differs subtly from \cite[Definition 7.9]{HJ22}, and the one given in the category of representations of quivers (see Remark~\ref{rem:difference_stalk}). This difference provides us an advantage in Section~6 when we characterize projective and injective functors in $\Add(\A, R\Mod)$ for more general categories $\A$ than the ones considered   in the literature.

The stalk functor is a generalization of  the $n$-sphere chain complex functor (see Example~\ref{ex:chain_stalk}) in the category of chain complexes, which is crucial in defining the classes $\tilde{\mathcal{F}}$ and $\tilde{\mathcal{G}}$. In Proposition~\ref{prop:adj2_stalk}, we show that the stalk functors induce  an adjoint triple $({ \rm c}, { \rm s}, { \rm k})$ between $\Add(\A, R\Mod)$ and $ \prod_{A \in \A} R_A \Mod$.  In Proposition~\ref{Prop:left_Derived_coeq}, we describe  ${\rm c}$-flat and ${ \rm k}$-coflat functors in $\Add(\A, R\Mod)$ in terms of  specific short exact sequences. This concrete characterization, as well as  providing us  a better understanding of ${\rm c}$-flat and ${ \rm k}$-coflat functors, is the key result in proving Theorem~\ref{them:proj-inj}. In the case of chain complexes, the short exact sequences given in Proposition~\ref{Prop:left_Derived_coeq}  correspond to the canonical short exact sequences
$${\rm c}_{i+1}(X) \longrightarrow X_i \longrightarrow {\rm c}_{i+1}(X) $$ 
$$ {\rm k}_{i}(X) \longrightarrow X_i \longrightarrow {\rm k}_{i-1}(X)$$
(for functors ${ \rm c}_i$ and ${ \rm k}_i$ on chain complexes, see Example~\ref{ex:chain_stalk}). Using Corollary~\ref{setup:3} and Proposition~\ref{Prop:left_Derived_coeq}, we obtain the following result.
\vspace{5mm}

\textbf{Theorem \ref{cor:cot_final}:}  If  $\prod_{A \in \A} \mcF_A$  and $\prod_{A \in \A} \mcG_A$ are a generating  and cogenerating classes in $\prod_{A \in \A} R_A \Mod$, respectively,  then
\begin{align*}
{ \rm s}( \prod_{A \in \A} \mcF_A)^\perp&=\biggl\{ X: \A \rightarrow R \Mod \biggm \vert \
\begin{array}{l}
{ \rm k}_A(X) \in \mcG_A,
X(A) \longrightarrow    \ker (\mu_{A,X})\\
\textrm{is an admissible epic in } R_A\Mod,\ \forall A \in \A
\end{array} \biggr \},\\
^\perp{ \rm s}( \prod_{A \in \A}  \mcG_A)&=\biggl \{ X: \A \rightarrow R \Mod \biggm \vert \
\begin{array}{l}
 { \rm c}_A(X) \in \mcF_A,  \coker (\phi_{A,X}) \hookrightarrow X(A)\\
\textrm{is an admissible monic in  } R_A\Mod,\ \forall A \in \A
\end{array}
\biggr \}.
\end{align*}

\vspace{4mm}

In Section~6, we focus on two trivial cotorsion pairs 
$$(\prod_{A \in \A} \mathcal{E}_A\mbox{-Proj}, \prod_{A \in \A} R_A \Mod), \quad \textrm{ and } \quad (\prod_{A \in \A} R_A \Mod, \prod_{A \in \A} \mathcal{E}_A\mbox{-Inj} )$$ to  provide an intrinsic caracterization of projective and injective objects in the exact category $(\Add(\A, R \Mod); \mathcal{E})$. 
\vspace{5mm}

 \textbf{Theorem~\ref{them:proj-inj}:}
\begin{enumerate}[(i)]
\item If Condition \ref{condition1} is satisfied, then  
$$\mathcal{E} \mbox{-Proj} = {}^\perp{ \rm s}( \prod_{A \in \A}  R_A \Mod).$$

\item If Condition \ref{condition2} is satisfied, then  
$$\mathcal{E} \mbox{-Inj} = { \rm s}( \prod_{A \in \A}  R_A \Mod)^\perp$$
\end{enumerate}
 It is worth pointing out   that  Theorem~\ref{them:proj-inj} provides a characterization of relative projective and injective functors. Moreover, in the absolute exact structure case, it applies to  more general categories $\A$ than the ones considered   in \cite[Theorem 7.29]{HJ22} (see Example~\ref{ex:non-left rooted}) and covers both a Morita context ring with zero trace ideals (as a particular case, formal triangular matrix rings) and  left/right rooted  categories. Therefore, our result extends the results in  \cite{HV99}, \cite{HV00} and \cite[Theorem1]{CT13}.

\section{Preliminares}
In this section, we review certain basic notions and results that  are essential for our study.
\subsection{Product category.} Let $\{ \C_i \}_{i \in I}$ be a family of additive categories. \textit{ The product category}, denoted by $\prod_{i\in I} \C_i$, is the category with objects

$$
\mbox{Ob}( {}\prod_{i\in I} \C_i ):= \prod_{i \in I} \mbox{Ob}( \C_i ),
$$
and morphisms
$$ \Hom_{  \prod_{i \in I} \C_i } ((C_i)_{i \in I}, (C'_i)_{i \in I}):=\prod_{i \in I} \Hom_{\C_i}(C_i,C_i').$$

For any $i \in I$, $\pi_i$ and $\iota_i$ denote the canonical $i$th  projection and inclusion functors, respectively,
$$\pi_i: \prod_{i \in I} \C_i \longrightarrow \C_i,$$
$$\iota_i: \C_i \longrightarrow \prod_{i \in I} \C_i. $$

\subsection{Exact categories.} An {\em exact structure on an additive category $ \mathcal{C}$}   is an isomorphism-closed collection $\mathcal{E}$ of distinguished kernel-cokernel pairs  $(i,p)$ in $\mathcal{C}$
$$ \mathbb{E} : \xymatrix {  A \ar@{^{(}->}[r]^i & B \ar@{->>}[r]^p & C},$$
which are called \textit{admissible monic} and \textit{admissible epic} of $\mathbb{E}$, respectively,
satisfying:
\begin{enumerate}[(i)]
\item For every object $A \in \mathcal{C},$  the identity morphism $\id_A$ is both an admissible monic and an admissible epic;
\item Admissible monics and admissible epics are closed under compositions;
\item The class $\mathcal{E}$ is closed under  pushout  and pullback   along a morphism in $\mathcal{C}$.
\end{enumerate}

$\mathcal{C}$ equipped with an exact structure $\mcE$, denoted by $(\mathcal{C}; \mathcal{E})$, is called an  \emph{exact category}. A kernel-cokernel pair  $(i,p)$, or $\mathbb{E}$, in $\mcE$ is  called a {\em short exact sequence } in $(\mcC; \mathcal{E})$.

A sequence in $\mathcal{C}$ of the form
$$ \mathbb{E}: \xymatrix{  A \ar@{^{(}->}[r] & B_n \ar[r] & \cdots \ar[r] & B_1 \ar@{->>}[r] & C}  $$
is said to be \textit{exact} if there exist short exact sequences $ \mathbb{E}_i: \xymatrix{  Z_{i} \ar@{^{(}->}[r]  & B_i \ar@{->>}[r] & Z_{i-1}}  $ in $(\mcC;\mcE)$, $ 1 \leq i \leq  n$, with $Z_0=C$ and $Z_n=A$ in such a way that the morphism ${B_{i+1} \longrightarrow B_{i}}$ is factorized as $\xymatrix{  B_{i+1}  \ar@{->>}[r] & Z_{i}   \ar@{^{(}->}[r]  & B_{i}}$ for every $1 \leq i \leq n-1.$ In other words, $\mathbb{E}=\mathbb{E}_n \circ \cdots \circ \mathbb{E}_1$. We refer to \cite{Buh} for a detailed treatement  on exact categories.

We recall the following result which will be used  in Section 4. 

\vspace{2mm}

\subsection{Proposition.}\label{prop:sum_admissible}\cite[Corollary 2.14]{Buh}. Let  $(\C; \mathcal{E})$ be an exact category. For any morphism  $f: A \rightarrow B$ in $\C$,  the following holds:
\begin{enumerate}[(i)]
\item  If $g: C \rightarrow B$ is an admissible epic in $(\C; \mathcal{E})$,  then so is the induced morphism $(f,g): A \oplus C \rightarrow B$.

\item  If $f: A \rightarrow C$ is an admissible monic in $(\C; \mathcal{E})$,  then so is the induced morphism $(f,g): A \rightarrow B \oplus C$.
\end{enumerate}

\subsection{Cotorsion pairs.}\label{subsection.Ext-orthogonality}  For any integer $n \geq 1$, the exact structure   $ \mcE$ on $\mcC$ leads to   the so-called \textit{$n$th Yoneda extension functor} $$\Ext^n_{\mathcal{E}} : \mcC^{\textrm{op}} \times \mcC \to \operatorname{Ab},$$ where
 $\Ext^n_{\mathcal{E}}(C,A)$ denotes the abelian group  of equivalence classes of  all exact sequences in $(\mcC;\mcE)$ of the form
$$ \mathbb{E}: \xymatrix{  A \ar@{^{(}->}[r] & B_n \ar[r] & \cdots \ar[r] & B_1 \ar@{->>}[r] & C}. $$

For given morphisms $i: C' \to C$,  $j: A \to A'$ in $\mcC$, the morphisms
$$\Ext^n_{\mathcal{E}}(i,A):\Ext^n_{\mathcal{E}}(C,A)\to\Ext^n_{\mathcal{E}}(C',A), \quad \mathbb{E} \leadsto \mathbb{E}i$$
$$\Ext^n_{\mathcal{E}}(C,j): \Ext^n_{\mathcal{E}}(C,A)\to \Ext^n_{\mathcal{E}}(C, A'), \quad \mathbb{E} \leadsto j\mathbb{E},$$
where  $\mathbb{E} i$ and $j \mathbb{E}$ denote the pullback and pushout of $\mathbb{E}$ along $i$ and $j$, respectively. So we have
$$\Ext^n_{\mathcal{E}}(i,j)= \Ext^n_{\mathcal{E}}(i,A') \circ \Ext^n_{\mathcal{E}}(C,j) \cong \Ext_{\mathcal{E}}^n(C',j) \circ \Ext_{\mathcal{E}}^n(i,A),$$
that is, $(j\mathbb{E})i \equiv j( \mathbb{E} i)$.

For a given class $\mathcal{S}$ of objects in $\mcC$, we let
\begin{align*}
\mathcal{S}^\perp & = \{ A \in \mcC \mid \Ext^1_{\mathcal{E}}(S,A)=0, \textrm{ for every } S \in \mathcal{S} \};\\
^\perp \mathcal{S} & = \{ A \in \mcC \mid \Ext^1_{\mathcal{E}}(A,S)=0,  \textrm{ for every } S \in \mathcal{S}\}.
\end{align*}

A pair $( \mathcal{F}, \mathcal{G})$ of classes of objects in $\mcC$ is said to be \textit{cotorsion pair} if $ \mathcal{F} ^\perp= \mathcal{G}$ and $\mathcal{F}= {}^\perp \mathcal{G}$. It is said to be \textit{hereditary} if $ \Ext^n_{\mathcal{E}}(F,G)=0$ for every $F \in \mathcal{F}$, $G\in \mathcal{G}$ and  $n \geq 1$.

For a  class $\mathcal{S}$ of objects in $\mcC$, the cotorsion pair $( \mathcal{F}, \mathcal{G})$ is said to be \textit{generated  by} $\mathcal{S}$ if $\mathcal{S}^\perp= \mathcal{G}$; is said to be  \textit{cogenerated by } $\mathcal{S}$ if $ \mathcal{F}={}^\perp \mathcal{S}$.

 The cotorsion pair $( \mathcal{F}, \mathcal{G})$ is said to \textit{have enough injectives} if for every object $A $ of $\mcC$, there exists a  short exact sequence in $(\mcC; \mathcal{E})$ of the form
$$
  \xymatrix{  A \ar@{^{(}->}[r] & G  \ar@{->>}[r] & F}
$$
with $G \in  \mathcal{G}$ and $F\in \mcF$; it is said to
\textit{have enough projectives} if for every object $A $ of $\mcC$, there exists a  short exact sequence in $(\mcC; \mathcal{E})$ of the form
$$
  \xymatrix{  G' \ar@{^{(}->}[r] & F'  \ar@{->>}[r] & A},
$$
where $G' \in  \mathcal{G}$ and $F' \in \mcF$; it is said to be
\textit{complete} if it has enough projectives and injectives.

\subsection{Projective and injective objects.} An object $P$ in $(\mcC; \mcE)$ is said to be \textit{projective} if for every short exact sequence  $\mathbb{E}$ in $\mcE$, the sequence $\Hom(P, \mathbb{E})$ of abelian groups is exact. The class of all projective objects in $(\mcC;\mcE)$ is denoted by $\mathcal{E}\mbox{-Proj}$.

$(\mcC; \mcE)$ is said to have \textit{enough projectives }if $(\mathcal{E}\mbox{-Proj}, \mcC )$ is a complete cotorsion pair, equivalently, for every object $A$ in $\mcC$ there exists an admissible epic $P \twoheadrightarrow A$ with  $P \in \mathcal{E}\mbox{-Proj}$.

Dually, an \textit{injective object} in  $(\mcC;\mcE)$ is defined, and the class of injective objects in $(\mcC;\mcE)$ is denoted by $\mathcal{E}\mbox{-Inj}$. $\mcC$ is said to have \textit{enough injectives} if for every object $A$ in $\mcC$ there exists an admissible monic $A \hookrightarrow E$ with  $E \in \mathcal{E}\mbox{-Inj}$.

\subsection{Left and right derived functors.}\label{subsec:left-right-derived}\cite[Section 10, Remark 12.12]{Buh} Let $\mathcal{F}$ be an extension closed class of objects in $(\mcC; \mcE)$ satisfying
\begin{enumerate}[(i)]
\item $\mcF$ is a \textit{generating class}, that is,  for every object $A$ in $\mathcal{C}$, there exists an admissible epic $F \twoheadrightarrow A$ with $F$ in $\mathcal{F}$;
\item $\mathcal{F}$ is closed under kernels of admissible epics, that is, if $\xymatrix{  F' \ar@{^{(}->}[r] & F  \ar@{->>}[r] & F''}$ is a short exact sequence in $(\mcC; \mcE) $ with objects $F,F''$ in $\mathcal{F}$, then so is $F'$.
\end{enumerate}

So for any object $A$ in $\mcC$ there exists an exact sequence in $(\mcC; \mcE)$ of the form
$$\xymatrix{\cdots \ar[r] &F_2 \ar[r] & F_1 \ar[r] & F_0 \ar@{->>}[r] & A }$$
with $F_i \in \mathcal{F}$ for every $i \geq 0$.
The deleted chain complex
$$F_*^A: \xymatrix{\cdots \ar[r] &F_2 \ar[r] & F_1 \ar[r] & F_0 \ar[r]& 0 }$$ is called  a $\mathcal{F}$-\textit{resolution}  of $A$. If an additive functor $\rm p: \mathcal{C} \longrightarrow \mathcal{D}$, where $\mathcal{D}$ is an abelian category, preserves short exact sequences in  $(\mcC; \mcE)$ with objects in $\mathcal{F}$, then for any  integer $i \geq 0 $,  the  so-called $i$th \textit{left derived functor  of $\rm p$ relative to $\mcF$}, denoted by $\rm L_i^{\mcF}p$, exists and is  the functor from $ \mathcal{C}$ to $\mathcal{D}$ defined by  $\rm{L_i^\mcF p}(A):= \operatorname{H}_i({ \rm p} (F^A_*))$, the  $i$th homology of the chain complex ${  \rm p}(F_*^A)$. It is well known that
 any short exact sequence in $(\mcC; \mcE)$
$$ \mathbb{E}:  \xymatrix {  A \ar@{^{(}->}[r]^i & B \ar@{->>}[r]^p & C}  $$
induces the following long exact sequence in $\mathcal{D}$
$$\scalebox{0.97}{\xymatrix{\cdots \ar[r] &  { \rm L_1^{\mcF} p}(B) \ar[r] & { \rm L_1^{\mcF} p}(C) \ar[r] & { \rm L_0^{\mcF} p}(A) \ar[r]& { \rm L_0^{\mcF} p}(B) \ar[r] & { \rm L_0^{\mcF} p}(C) \ar[r] & 0}} $$
and if ${ \rm p}$ is right exact, then ${ \rm L_0}^{\mcF}  { \rm p} \cong \rm p$.

If $\mcF$ satisfies dual properties, then every object in $\mathcal{C}$ has a $\mathcal{F}$-\textit{coresolution}, and the $i$th \textit{right derived functor}, denoted by ${ \rm R^i_{\mcF}p}: \mcC \longrightarrow \mathcal{D}$, of $ { \rm p}$  exists and is defined dually.

\subsection{} Returning to the product categories, for a  family $\{ (\C_i;\mathcal{E}_i) \}_{i \in I}$  of exact  categories, \textit{the product exact structure} $\prod_{i \in I} \mathcal{E}_i$ on  the product category $\prod_{i\in I} \C_i$  consists of $I$-tuples $(\mathbb{E}_i)_{i\in I}$ of short exact sequences $\mathbb{E}_i$ in $(\mathcal{C}_i; \mcE_i)$.
So we have
\begin{equation}\label{ext-product-category}
 \Ext^n_{  \prod_{i \in I} \mcE_i } ((C_i)_{i \in I}, (A_i)_{i \in I}):=\prod_{i \in I} \Ext^n_{\mcE_i}(C_i,A_i).
\end{equation}
The following  is an easy observation, but for the reader's convenience, we provide its proof.

\subsection{Lemma.}\label{lemma:product_cot_pair} Any (hereditary) cotorsion pair in $(\prod_{i \in I} \C_i ; \prod_{i \in I} \mathcal{E}_i)$ is of the form
$$(\  {} \prod_{i \in I} \mcF_i,  \prod_{i \in I} \mcG_i),$$
where  $(\mcF_i, \mcG_i)$ is a (hereditary) cotorsion pair in $(\mathcal{C}_i; \mathcal{E}_i)$ for every $i \in I$.

\begin{proof}
Let $\{ (\mcF_i, \mcG_i) \}_{i \in I}$ be a family of   (hereditary) cotorsion pairs in $(\mathcal{C}_i; \mathcal{E}_i)$. From \eqref{ext-product-category}, we already  have
$ \prod_{i \in I} \mcF_i \subseteq {} ^\perp ( \prod_{i \in I} \mcG_i)$ and
$ \prod_{i \in I} \mcG_i \subseteq\  ( \prod_{i \in I} \mcF_i)^\perp$. Let  $ ( C_i)_{i \in I}  $ be an object   in
$ ^\perp (\prod_{i \in I} \mcG_i)$.  It is clear that for any object  $G_{i_0} $  in $ \mcG_i$, ${\iota_{i_0}(  G_{i_0} ) \in \prod_{i \in I } \mcG_{i}}$, and by assumption,
 we have
$$0=
 \Ext^1_{\prod_{i\in I} \mcE_i}(\ (C_i)_{i\in I} , \iota_{i_0}(G_{\iota_0}))=\Ext^1_{\mcE_{i_0}}(C_{i_0}, G_{i_0}).
	$$
Therefore, for every $i \in I$, $C_{i} \in\ ^\perp \mcG_{i}=\mcF_{i}.$ Similarly, we have  ${( \prod_{i \in I} \mcF_i )^\perp \subseteq \   \prod_{i \in I} \mcG_i}$.

Now, suppose that $(\mcF, \mcG)$ is a cotorsion pair in $(\prod_{i \in I} \C_i ; \prod_{i \in I} \mathcal{E}_i)$. For every $i \in I$, we let $\mcF_i:=\pi_i(\mcF)$ and $\mcG_i:=\pi_i(\mcG)$. From \eqref{ext-product-category}, it is easy to observe that $\mcF_i \subseteq {}^\perp \mcG_i$ and $\mcG_i \subseteq \mcF_i^\perp $ for every $i \in I$. Now, let $C_{i_0} \in   {}^\perp \mcG_{i_0}$ for some $i_0 \in I$. For every $(G_i)_{i\in I} \in \mathcal{G}$, we have  $G_{i_0} \in \mathcal{G}_{i_0}$ and
$$\Ext^1_{\prod_{i\in I} \mcE_i}(\iota_{i_0} (C_{i_0}), (G_i)_{i\in I})=\Ext^1_{\mcE_{i_0}}(C_{i_0}, G_{i_0})=0,$$
therefore,  $\iota_{i_0}(C_{i_0}) \in \mcF$. It implies that $C_{i_0}=\pi_{i_0}( \iota_{i_0} (C_{i_0})) \in \mathcal{F}_{i_0}$. Similarly, we have $\mathcal{G}_i=\mathcal{F}_i^\perp$ for every $i \in I$.

The hereditary condition follows easily.
\end{proof}

\section{$\Ext$ isomorphisms.}
 Throughout, $\mathcal{C}$ and $\mathcal{D}$ denote weakly idempotent complete exact categories equipped with  exact structures $\mathcal{E}$ and $\mathcal{E}'$, respectively. So we omit the subscript in the corresponding $n$th Yoneda extension functors. We also fix  an adjoint pair $({ \rm q}, { \rm t})$ of functors between  $\mathcal{C}$ and $\mathcal{D}$
$$
\xymatrix{\mathcal{C} \ar@/^1em/[rrr]^{ \rm q} &&& \mathcal{D} \ar@/^1em/[lll]^{ \rm t} }.
$$

We let $\eta: \id_{\mcC} \Rightarrow { \rm tq}$ and $\xi: { \rm qt } \Rightarrow \id_{\mathcal{D}}$ denote the unit and the counit of the adjunction $({ \rm q}, { \rm t})$. It is  well known that ${ \rm q}$ preserves all  colimits   while ${ \rm t}$ preserves all limits.

In this section, we investigate sufficient and necessary conditions for which the functors ${ \rm q}$ and ${ \rm t}$  interact together inside the functor $\Ext^n$ in the same way as they do inside of the $\Hom$ functor.
The following presents  a more general form of the result \cite[Lemma 5.1]{HJ19}.

\subsection{Lemma.}\label{lemma:ext_adjoint_pair}
 For any objects $C$ and $D$ in $\mathcal{C}$ and $\mathcal{D}$, respectively, and an integer $n \geq 1$, the following hold:
\begin{enumerate}[(i)]
	\item If  ${ \rm t}$ preserves exact sequences in $\Ext^n({ \rm q}(C),D)$, then there exists a canonical morphism
$$\Theta_{CD}^n:\ \Ext^n({ \rm q}(C), D) \longrightarrow\Ext^n(C, { \rm t}(D)).$$
If $n=1$, then $\Theta_{CD}^1$ is a monomorphism.

\item If  ${ \rm q}$ preserves exact sequences in $\Ext^n(C, { \rm t}(D))$, then there exists a canonical morphism
$$\Omega_{CD}^n:\ \Ext^n(C, { \rm t}(D)) \longrightarrow \Ext^n({ \rm q}(C), D).$$
If $n=1$, then $ \Omega_{CD}^1$ is a monomorphism.

\item If functors ${ \rm q}$ and ${ \rm t}$ satisfy conditions given in (i)-(ii), then
$$\Theta_{CD}^n:\ \Ext^n({ \rm q}(C), D) \longrightarrow\Ext^n(C, { \rm t}(D))$$
is an isomorphism with the inverse  $\Omega_{CD}^n$.

\end{enumerate}

\begin{proof}
\begin{enumerate}[(i)]
	\item	
Consider an exact sequence   $\mathbb{E}$ in $\Ext^n({ \rm q}(C), D)$. By assumption, ${ \rm t}(\mathbb{E})$ is an  exact sequence  in  $\Ext^n(tq(C), { \rm t}(D))$.
Taking the pullback of ${ \rm t}(\mathbb{E})$ along  the unit morphism $\eta_C:\ C \longrightarrow tq(X)$, we obtain  an  exact sequence   ${ \rm t}(\mathbb{E})\eta_C $ in $\Ext^n(C, { \rm t}(D))$ together with the pullback diagram  ${ \rm t}(\mathbb{E}) \eta_C \rightarrow { \rm t}(\mathbb{E})$.

Since taking pullback and ${ \rm t}$ are additive functors, the  assignment
$$\Theta^n_{CD}:\ \Ext^1({ \rm q}(C), D) \longrightarrow \Ext^n(C, { \rm t}(D)), \quad \mathbb{E} \rightsquigarrow { \rm t}(\mathbb{E})\eta_C$$
is well-defined.

\sloppy Now, we show that $\Theta^1_{CD}$ is a monomorphism. Suppose that for some short exact sequence $\mathbb{E}$ in $\Ext^1({ \rm q}(C), D)$, the short exact sequence ${ \rm t}(\mathbb{E})\eta_C$  is trivial in $ \Ext^1(C, { \rm t}(D))$. As ${ \rm q}$ is additive,   ${ \rm q}({ \rm t} (\mathbb{E}) \eta_C)$ is trivial, as well. On the other hand, we consider the composition of the following canonical morphisms  between  sequences
	\begin{equation}\label{composit}
	{ \rm q}({ \rm t}(\mathbb{E}) \eta_C) \longrightarrow qt(\mathbb{E}) \longrightarrow \mathbb{E},
	\end{equation}
	where $qt(\mathbb{E}) \longrightarrow \mathbb{E}$ is the morphism induced from  the counit natural transformation $\xi$. The last column of the composition given in \eqref{composit} is the  composition $$\xi_{{ \rm q}(C)} \circ { \rm q}(\eta_C): { \rm q}(C) \longrightarrow qtq(C) \longrightarrow { \rm q}(C),$$ which  is identity. Hence, the short exact sequence $\mathbb{E}$ is in fact     a pushout of the trivial short exact sequence ${ \rm q}({ \rm t}(\mathbb{E})\eta_C) $. As a consequence, $\mathbb{E}$ is a trivial short exact sequence.

\item Similar to the proof of the first statement.

\item From (i) and (ii), we have assignments $\Theta_{CD}^n$ and $\Omega_{CD}^n$. We claim that
$$\Omega_{CD}^n \circ \Theta_{CD}^n =\id \quad \textrm{ and } \quad  \Theta_{CD}^n \circ \Omega_{CD}^n=\id. $$
 By duality of arguments, we only show the equality $\Theta_{CD}^n \circ \Omega_{CD}^n=\id $.

 Let $\mathbb{E}$ be any  exact sequence in $\Ext^n({ \rm q}(C),D)$.  By construction, we have $\Theta^n_{CD}(\mathbb{E})={ \rm t}(\mathbb{E}) \eta_C$ together with the pullback diagram $ { \rm t}(\mathbb{E}) \eta_C \longrightarrow { \rm t}(\mathbb{E})$. Applying ${ \rm q}$ to the aforementioned pullback diagram, we still have a commutative  diagram ${{ \rm q}({ \rm t}(\mathbb{E}) \eta_C) \longrightarrow qt(\mathbb{E}) }$ with  the rightest vertical arrow ${ \rm q}({ \rm t}( \eta_C)) $. By assumption, ${ \rm q}({ \rm t}(\mathbb{E}) \eta_C) $ is an  exact sequence in $\mathcal{D}$. Composing it with the counit morphism, we have a commutative diagram with canonical morphisms
	\begin{equation*}
	{ \rm q}({ \rm t}(\mathbb{E}) \eta_X) \longrightarrow qt(\mathbb{E}) \longrightarrow \mathbb{E},
	\end{equation*}
	The last column of the composition  is the  composition $$\xi_{p(C)} \circ { \rm q}(\eta_C): { \rm q}(C) \longrightarrow qtq(C) \longrightarrow { \rm q}(C),$$ which  is identity. And the left vertical map is $\xi_D$. By universal property of pushouts, the commutative diagram $	{ \rm q}({ \rm t}(\mathbb{E}) \eta_X)  \longrightarrow \mathbb{E}$ is factorized as
	$$ { \rm q}({ \rm t}(\mathbb{E}) \eta_X) \longrightarrow \xi_D { \rm q}({ \rm t}(\mathbb{E}) \eta_X) \longrightarrow \mathbb{E}$$
where both $\Omega^n_{CD} \Theta^n_{CD}( \mathbb{E})=\xi_D { \rm q}({ \rm t}(\mathbb{E}) \eta_X) $	and $\mathbb{E}$ belongs to $\Ext^n({ \rm q}(C),D)$. So $\Omega^n_{CD} \Theta^n_{CD}( \mathbb{E}) \equiv \mathbb{E}$.

\end{enumerate}

\end{proof}

As an aside, it is easy to verify that if $1 \leq m \leq n$ and ${ \rm t}$ preserves exact sequences in $\Ext^n({ \rm q}(C),D)$, then it  preserves also exact sequences in $\Ext^m({ \rm q}(C),D)$, and therefore,  there exists a morphism $\Theta_{CD}^m:\ \Ext^m({ \rm q}(C), D) \longrightarrow\Ext^m(C, { \rm t}(D)).$ The dual affirmation is valid for ${ \rm q}$.

As pullback and pushout diagrams satisfy the universal property, and the (co)unit of an  adjunction is a natural transformation, we have the following result.
\subsection{Corollary.}\label{corol:ext3}
Let $n \geq 1$. We have the following:
\begin{enumerate}[(i)]
	\item If  ${ \rm t}$ preserves  exact sequences in $\Ext^n({ \rm q}(-),-)$, then $\Theta^n:=\{ \Theta_{CD}^n\}_{C\in \mathcal{C}, D\in \mathcal{D} }$ is a natural transformation
$$\Theta^n:\ \Ext^n({ \rm q}(-), -) \Longrightarrow\Ext^n(-, { \rm t}(-)).$$

\item If  ${ \rm q}$ preserves  exact sequences in $\Ext^n(-, { \rm t}(-))$, then $\Omega^n:=\{ \Omega^n_{CD}\}_{C\in \mathcal{C}, D\in \mathcal{D}}$ is a natural transformation
$$\Omega^n:\ \Ext^n(-, { \rm t}(-)) \Longrightarrow \Ext^n({ \rm q}(-), -).$$

\item If functors ${ \rm q}$ and ${ \rm t}$ satisfy conditions given in (i)-(ii), then the natural transformation
$$\Theta^n:\ \Ext^n({ \rm q}(-), -) \Longrightarrow\Ext^n(-, { \rm t}(-))$$
is an isomorphism with the inverse  $\Omega^n$.

\end{enumerate}

\subsection{Remark} The condition that ${ \rm t}$ preserves  exact sequences in $\Ext^n({ \rm q}(-),-)$ for some $n \geq 1$ is equivalent to the exactness of  ${ \rm t}$; see Lemma \ref{lemma:lem}. Therefore, if   both functors ${ \rm t}$ and ${ \rm q}$ are exact, then the natural transformation $$\Theta^n:\ \Ext^n({ \rm q}(-), -) \Longrightarrow\Ext^n(-, { \rm t}(-))$$
is an isomorphism for every $n \geq 1$.

\subsection{Example.}\label{ex:chain_complex1}
We let $(\mathbf{C}(\mathcal{C});\mathbf{C}(\mathcal{E}))$ denote the category of chain complexes of objects in  $\mathcal{C}$ equipped with the degrewise exact structure. For every $i\in \mathbb{Z}$, there is an adjoint triple $( \operatorname{D}^i, \operatorname{ev}_i, \operatorname{D}^{i+1})$
	$$
\xymatrix{\mathcal{C} \ar@/^2em/[rrr]^{\operatorname{D}^{i}} \ar@/_2em/[rrr]_{\operatorname{D}^{i+1}} &&& \mathbf{C}(\mathcal{C}) \ar[lll]^{\operatorname{ev}_i} },
$$
where $\operatorname{ev}_i $ denotes the $i$th evaluation functor, and $\operatorname{D}^i(C)$ is the chain complex with $C$ at $i$th and $(i-1)$th positions, and $0$ elsewhere. For every $i\in \mathbb{Z}$, the functors $\operatorname{D}^i$ and $\operatorname{ev}_i$ are known to be exact, and therefore, we have  natural isomorphism
$$\Ext^n( \operatorname{D}^i(-), -) \cong \Ext^n(-, \operatorname{ev}_i (-)) \quad \textrm{ and } \quad   \Ext^n( \operatorname{ev}_i (-),-) \cong   \Ext^n( -,\operatorname{D}^{i+1}(-)) $$
for every $n \geq 1$.

\subsection{Example.} \sloppy Let $X$ be a semi-separated scheme. For any affine subscheme ${\iota:U \hookrightarrow X}$,  we have the adjoint pair $(\operatorname{res}, \iota_*)$ of functors

$$
\xymatrix{   \mathfrak{Qco}(X)  \ar@/^1em/[rrr]^{ \operatorname{res}}   &&& \mathfrak{Qco}(U)  \ar@/^1em/[lll]^{ \iota_*}      }
$$
between categories of quasi-coherent sheaves on $X$ and $U$. 
Both the restriction $\operatorname{res}$ and the direct image functor $\iota_*$  are exact, and therefore, for every $n \geq 1$
$$\Ext^n(\operatorname{res}(-),-) \cong \Ext^n(-, \iota_*(-)).$$

\vspace{3mm}

From now on, even though we  present statements together with their dual versions, we only prove the first statement as their arguments can be dualized easily. We also use notations presented  in Lemma \eqref{lemma:ext_adjoint_pair}. 

 Now we investigate conditions for  which $\Theta_C^n$ and $\Omega_D^n$ are isomorphisms for fixed objects $C \in \mathcal{C}$ and $D \in \mathcal{D}$. As a first step, we study the case $n=1$. For it, we introduce the following definition.

\subsection{Definition.}\label{def:flat-coflat} We call an object $C$ in $\mathcal{C}$ $\rm q$-\textit{flat} if every short exact sequence in $\mathcal{C}$ ending with $C$ remains exact under $ \rm q$. Dually, a $\rm t$-\textit{coflat} object in $\mathcal{D}$ is defined.

\subsection{Remark.} In \cite{Pre10}, for given two proper classes $\mathcal{P}$ and $\mathcal{R}$ in an abelian category $\mathcal{B}$,  the notions of $\mathcal{P}$-$\mathcal{R}$-flat and divisible object have been introduced. Namely,  an object $B$ in $\mathcal{B}$ is called $ \mathcal{P}$-$\mathcal{R}$-flat
 if every short exact sequence
in $\mathcal{P}$ ending at $B$ belongs to $\mathcal{R}$. $\mathcal{P}$-$\mathcal{R}$ divisibles are defined dually. In our case, since $ { \rm q}$ is a left adjoint functor, it can be easily proved that the class $\overline{ \mathcal{E}}$ consisting of short exact sequences in $\mathcal{E}$ which remains exact under ${ \rm q}$ is an exact substructure of $\mathcal{E}$. So an object $C$ in $\mathcal{C}$ is ${ \rm q}$-flat if and only if it  is $\mathcal{E}$-$\overline{\mathcal{E}}$-flat in the sense of  \cite{Pre10}. The dual statment holds for ${p}$-coflat and divisible objects.

\vspace{3mm}

The following lemma will be used for proving Proposition \ref{prop:adjoints}. 

\subsection{Lemma.}\label{lemma:lem}
 For given objects $C \in \mcC$ and $D \in \mathcal{D}$, the following hold:

	\begin{enumerate}[(i)]
	\item  The functor ${ \rm q}$ preserves short exact sequences in $\Ext^1(C,{ \rm t}(-))$ if and only if  $C$ is $\rm q$-\textit{flat}. 

\item  The functor ${ \rm t}$ preserves short exact sequences in $\Ext^1({ \rm q}(-),D)$  if and only if   $D$ is $\rm t$-\textit{coflat}. 
\end{enumerate}

\begin{proof}
The ``if'' condition of the statement (i) is immediate. For the ''only if'' part, suppose that ${ \rm q}$ preserves short exact sequences in $\Ext^1(C,{ \rm t}(-))$. Consider a short exact sequence in $\mathcal{C}$
$$ \mathbb{E}: \xymatrix{ A \ar@{^{(}->}[r]^i & B \ar@{->>}[r]^p & C}.$$
By taking pushout of $\mathbb{E}$ along the unit morphism $\eta_A$, we have a short exact sequence in $\mathcal{C}$
$$ \eta_A \mathbb{E}: \xymatrix{ { \rm tq}( A) \ar@{^{(}->}[r] & B' \ar@{->>}[r]& C}$$
together with the pushout diagram $\mathbb{E} \longrightarrow \eta_A \mathbb{E}$. By assumption, ${ \rm q}(\eta_A \mathbb{E}) $ is a short exact sequence in $\mathcal{D}$. Taking pushout of ${ \rm q}(\eta_A \mathbb{E}) $ along the morphism $\xi_{ { \rm q}(A)}$, we obtain a short exact sequence $\xi_{ { \rm q}(A)}{ \rm q}(\eta_A \mathbb{E}) $ together with the commutative diagram
$${ \rm q}(\mathbb{E}) \longrightarrow { \rm q}(\eta_A \mathbb{E}) \rightarrow \xi_{ { \rm q}(A)}{ \rm q}(\eta_A \mathbb{E})$$
$$
\xymatrix{
{ \rm q}(\mathbb{E}): \ar[d] &  { \rm q}(A) \ar[r]^{{ \rm q}(i)}  \ar@{=}[d] & { \rm q}(B) \ar[d] \ar[r]^{{ \rm q}(p)} & { \rm q}(C) \ar@{=}[d]\\
\xi_{ { \rm q}(A)}{ \rm q}(\eta_A \mathbb{E}):  &{ \rm q}(A) \ar@{^{(}->}[r] & X \ar@{->>}[r] &{ \rm q}(C)
}
$$

Since ${ \rm q}$ preserves cokernels, ${ \rm q}(p)$ is the cokernel of ${ \rm q}(i)$, and by \cite[Proposition~2.1.6]{Buh}, ${ \rm q}(i)$ is an admissible monic in $\mathcal{D}$, and therefore, ${ \rm q}(\mathbb{E})$ is a short exact sequence in $\mathcal{D}$.

\end{proof}

The following lemma,  which will be used to prove Theorem \ref{Prop:degrewise-cot.pair}, illustrates that (co)generating classes are sufficient for showing  (co)flatness of an object. 
\subsection{Lemma.}\label{lemma:exactt} Let $\mathcal{F}$ and $\mathcal{G}$ be classes of objects in $\mathcal{C}$ and $\mathcal{D}$, respectively. The followings hold:
\begin{enumerate}[(i)]
\item  If $ \mathcal{G}$ is a cogenerating class in $\mathcal{D}$, and  ${ \rm q}$ preserves short exact sequences in $\Ext^1(C,{ \rm t}(\mcG))$,  then   $C$ is a ${ \rm q}$-flat object.

\item \sloppy If    $\mathcal{F}$ is a generating class in $\mathcal{C}$, and   ${ \rm t}$  preserves short exact sequences in $\Ext^1({ \rm q}(\mcF),D)$,  then $D$ is a $\rm t$-coflat object.

\end{enumerate}	
\begin{proof}
 Let $\mathbb{E}: \xymatrix{  C'' \ar@{^{(}->}[r] & C'  \ar@{->>} [r] & C} $ be a short exact sequence in $\mathcal{C}$. By assumption, there exists an admissible monic $f:  { \rm q}(C'') \hookrightarrow G$ in $\mathcal{D}$ with $G \in \mcG$. Let $\overline{f}:=  {\rm t}(f) \circ \eta_{C''}: C'' \rightarrow  { \rm t}(G)$.
Consider the  pushout diagram $\mathbb{E}  \longrightarrow \overline{f} \mathbb{E}$. $\overline{f} \mathbb{E} $ is a short exact sequence in $\Ext^1(C, { \rm t}( G ))$, and  by assumption, ${ \rm q}( \overline{f} \mathbb{E} )$ is a short exact sequence in $\mathcal{D}$. Again, taking pushout of the short exact sequence ${ \rm q}(\overline{f} \mathbb{E}  )$ along the counit morphism $\xi_G :  { \rm qt}(G) \rightarrow G$, we obtain a short exact sequence $\xi_G { \rm q}(\overline{f} \mathbb{E}  )$. Consider the natural composition of sequences
$${ \rm q}( \mathbb{E})  \longrightarrow { \rm q}(  \overline{f}  \mathbb{E} ) \longrightarrow \xi_G { \rm q}( \overline{f}  \mathbb{E}  ).  $$
\sloppy Note that the morphism on the most left column of the above composite is $\xi_G \circ  {{ \rm q} ( \overline{f}) =f}$, which is an admissible monic in $\mathcal{D}$. Since ${ \rm q}$ preserves cokernels, the morphism ${{ \rm q}(C') \rightarrow {\rm q}(C)}$ is the cokernel of the morphism  ${ \rm q}(C'') \rightarrow {\rm q}(C')$, and by \cite[Proposition 2.16]{Buh},   ${ \rm q}(C'') \rightarrow {\rm q}(C')$ is an admissible monic in $\mathcal{D}$, and therefore, ${ \rm q} (\mathbb{E})$ is a short exact sequence in $\mathcal{D}$.
\end{proof}

\subsection{Proposition.}\label{prop:adjoints}
 For given objects $C \in \mcC$ and $D \in \mathcal{D}$, the following hold:
\begin{enumerate}[(i)]
	\item Suppose that  ${ \rm t}$ preserves short exact sequences in $\Ext^1({ \rm q}(C),-)$. Then the following are equivalent:
	
	\begin{enumerate}[(a)]
	\item The natural transformation $$\Theta_C^1:\ \Ext^1({ \rm q}(C), -) \Longrightarrow\Ext^1(C, { \rm t}(-)),$$
is an isomorphism.
	\item   $C$ is a ${ \rm q}$-flat object in  $\mcC$.
\end{enumerate}
\item Suppose that  ${ \rm q}$ preserves short exact sequences in $\Ext^1(-, { \rm t}(D))$. Then the following are equivalent:
	
	\begin{enumerate}[(a)]
	\item The natural transformation $$\Omega_D^1:\ \Ext^1(-, { \rm t}(D)) \Longrightarrow\Ext^1({ \rm q}(-),D),$$
is an isomorphism.
	\item  $D$ is a ${ \rm t}$-coflat object in $\mathcal{D}$.
\end{enumerate}
\end{enumerate}
\begin{proof}
The implication (b)$\Rightarrow$(a) in (i) follows from Lemma  \ref{lemma:ext_adjoint_pair}-(iii) and Lemma \ref{lemma:lem}-(i).
Now suppose that the natural transformation $ \Theta_C^1$ is an isomorphism. Consider a short exact sequence in $\mathcal{C}$
$$ \mathbb{E}: \xymatrix{ A \ar@{^{(}->}[r]^i & B \ar@{->>}[r]^p & C}.$$
By taking pushout of $\mathbb{E}$ along the unit morphism $\eta_A$, we have a short exact sequence in $\mathcal{C}$
$$ \eta_A \mathbb{E}: \xymatrix{ { \rm tq}( A) \ar@{^{(}->}[r] & B' \ar@{->>}[r]& C}$$
together with the pushout diagram $\mathbb{E} \longrightarrow \eta_A \mathbb{E}$. By assumption, there exists a short exact sequence $\mathbb{E}': \xymatrix{ { \rm q}(A) \ar@{^{(}->}[r] & X \ar@{->>}[r] & { \rm q} (C)} $  in $\mathcal{D}$  such that ${ \rm t}(\mathbb{E}') \eta_C \equiv \eta_A \mathbb{E}$. Note that, by assumption, ${ \rm t}(\mathbb{E}')$ is a short exact sequence in $\mathcal{C}$. Applying ${ \rm q}$, we have a commutative diagram  ${ \rm q}(\mathbb{E}) \longrightarrow { \rm q}(\eta_A \mathbb{E}) \equiv  { \rm q}({ \rm t}(\mathbb{E}') \eta_C )\rightarrow qt(\mathbb{E}').$
Composing it with the canonical diagram $qt(\mathbb{E}') \longrightarrow \mathbb{E}'$ induced from the counit natural transformation $\xi$, we have a  commutative diagram
$$
\xymatrix{
{ \rm q}(\mathbb{E}): \ar[d] &  { \rm q}(A) \ar[r]^{{ \rm q}(i)}  \ar@{=}[d] & { \rm q}(B) \ar[d] \ar[r]^{{ \rm q}(p)} & { \rm q}(C) \ar@{=}[d]\\
\mathbb{E}':  &{ \rm q}(A) \ar@{^{(}->}[r] & X \ar@{->>}[r] &{ \rm q}(C)
}
$$
Since ${ \rm q}$ preserves cokernels, ${ \rm q}(p)$ is the cokernel of ${ \rm q}(i)$, and by \cite[Proposition~2.1.6]{Buh}, ${ \rm q}(i)$ is an admissible monic in $\mathcal{D}$, and therefore, ${ \rm q}(\mathbb{E})$ is a short exact sequence in $\mathcal{D}$.
\end{proof}

	
The previous result  covers the one proven in \cite[Theorem 4.2]{AP22} for ${ \bf I}$-directed colimit functor, as shown in the following example.

\subsection{Example.}\label{ex:colimit}
Let $\mbox{Fun}( { \bf I},\mathcal{C})$ be the category of functors from ${ \bf I}$ to   $\mathcal{C}$ equipped with the degrewise exact structure $\mbox{Fun}( { \bf I}, \mathcal{E})$. We let $\kappa_{ \bf I} :\ \mathcal{C} \longrightarrow \mbox{Fun}( { \bf I},\mathcal{C})$ denote  the constant functor, that is, for a given object $C$ in $\mathcal{C}$,  $\kappa_{ \bf I}(C)(i):=C$ and $\kappa_{ \bf I}(C)(f):=\mbox{id}_C$ for every object $i$ and a morphism $f$ in ${ \bf I}$. It is a well-known fact that $\kappa_{ \bf I}$ is an exact functor, and  $\kappa_{ \bf I}$ has a left adjoint if and only if $\mathcal{C}$ has ${ \bf I}$-directed colimits.  We assume that  $\mathcal{C}$ has ${ \bf I}$-directed colimits. For a given functor  $F:{ \bf I} \rightarrow \mathcal{C} $, $\mbox{colim}_{\bf I}(F)$ denotes  its colimit. So the canonical natural transformation
$$\Ext^1(\mbox{colim}_{\bf I}(F), -) \Longrightarrow\Ext^1(F, \kappa_{ \bf I}(-))$$ is an isomorphism if and only if for every short exact sequence $\mathbb{E}$ in $\mbox{Fun}( { \bf I},\mathcal{C})$ ending with $F$, $\mbox{colim}_{\bf I}( \mathbb{E})$ is a short exact sequence in $\mathcal{C}$, as well. In particular, if $\mathcal{C}$ is an abelian category with ${ \bf I}$-directed colimits, then the functor of ${ \bf I}$-directed colimits in $\mathcal{C}$ is exact if and only if $\Ext^1(\mbox{colim}_{\bf I}(-), -) \Longrightarrow\Ext^1(-, \kappa_{ \bf I}(-))$ is an isomorphism; see \cite[Theorem 4.2]{AP22}.

As a particular case, if $\bf I$ is a discrete category, then $\mbox{Fun}( { \bf I},\mathcal{C}) \cong  \prod_{i \in { \bf I}} \mathcal{C}$, and $\mbox{colim}_{\bf I}(-) \cong \oplus_{i \in {\bf I}}$. So $\mathcal{C}$ has ${\bf I}$-indexed coproducts if and only if $\Ext^1(\oplus_{\bf I}(-), -) \cong \prod_{i \in {\bf I}}\Ext^1(-, -)$.

\vspace{3mm}

To obtain similar results for  $\Theta^n_C $  and $\Omega^n_D $ for any $n \geq 1$ in such generality, we would need to consider some lengthy technical conditions, such as being all syzygies of $C$ ${ \rm q}$-flat. Instead,   we prefer to use left and right derived functors (see Proposition \ref{lemma:derived_zero}).  For it, we prove the following two lemmas.

\subsection{Lemma.}\label{prop:heredit1} Let $\mathcal{F}$ and $\mathcal{G}$ be extension closed-classes of objects in $\mathcal{C}$ and $\mathcal{D}$ satisfying conditions (i)-(ii) and their duals given in \eqref{subsec:left-right-derived},  respectively. The followings hold:
\begin{enumerate}[(i)]
\item Suppose that  ${ \rm t}$  preserves short exact sequences in $\Ext^1({ \rm q}(\mcF),-)$.
The functor ${ \rm q}$ restricted on $\mcF$ is exact if and only if $\Theta_{F}^n$ is an isomorphism for every $F \in \mathcal{F}$ and $n \geq 1$.

\item  Suppose that  ${ \rm q}$ preserves short exact sequences in $\Ext^1(-,{ \rm t}(\mcG))$.
The functor  ${ \rm t}$ restricted on ${\mcG}$  is exact if and only if $\Omega_D^n$ is an isomorphism for every $G \in \mathcal{G}$ and $n \geq 1$.
\end{enumerate}	
\begin{proof}
Firstly,  note that the functor ${ \rm t}$ is an exact functor by Lemma \ref{lemma:exactt}. So the natural transformation $\Theta_{F}^n: \Ext^n({ \rm q}(F),-) \Longrightarrow \Ext^n(F,{ \rm t}(-))$ exists for every $F \in \mathcal{F}$ and $n \geq 1$.

The ''if'' statement now follows from Proposition \ref{prop:adjoints}. For the ''only if'' part,  suppose that ${ \rm q}$ preserves short exact sequences in $\mcC$ with objects in $\mathcal{F}$. We firstly prove  that ${ \rm q}$ preserves all short exact sequences in $\mathcal C$ ending with an object in $\mathcal{F}$.   Let
$$\mathbb{E}:\ \xymatrix{A \ar@{^{(}->}[r] &B \ar@{->>}[r] & F}$$
be a short exact sequence in $\mathcal{C}$ with an object $F$ in $\mathcal{F}$. By assumption, there exists an admissible epic $F_0 \twoheadrightarrow B$ with $F_0 \in \mathcal{F}$, which induces the following  commutative diagram
$$
 \xymatrix{\mathbb{E}': &
 F_1 \ar@{^{(}->}[r] \ar[d]_f & F_0 \ar@{->>}[d]\ar@{->>}[r] & F\ar@{=}[d]\\
\mathbb{E}: &  A \ar@{^{(}->}[r] &B \ar@{->>}[r] & F}
$$
with exact rows, and therefore, $f\mathbb{E}'=\mathbb{E}$. Again by assumption, $F_1 \in \mathcal{F}$, and hence, ${ \rm q}(\mathbb{E}')$ is a short exact sequence in $\mathcal{D}$. Since ${ \rm q}$ preserves all colimits, it preserves cokernels and pushout diagrams. So we have
$${ \rm q}(\mathbb{E})\cong { \rm q}(f \mathbb{E}') \cong { \rm q}(f) { \rm q}(\mathbb{E}'),$$ that is, ${ \rm q}(\mathbb{E})$ is a pushout of ${ \rm q}(\mathbb{E}')$ along ${ \rm q}(f)$.  It implies that ${ \rm q}(\mathbb{E})$ is a short exact sequence in $\mathcal{D}$, and by  Proposition \ref{prop:adjoints}, $\Theta_{F}^1$ is an isomorphism for every $F$ in $\mathcal{F}$.

Now, we prove the statement for $n=2$. The natural transformation
$$\Theta_F^2:\ \Ext^2({ \rm q}(F), -) \Longrightarrow\Ext^2(F, { \rm t}(-))$$ is an epimorphism
for every $F \in \mathcal{F}$. Indeed, by a similar argument just as done above, any  exact sequence  in $\Ext^2(F, { \rm t}(D))$ is equivalent to an exact sequence of the form
\begin{equation}\label{exactF}
\mathbb{E}:\ \xymatrix{ { \rm t}(D) \ar@{^{(}->}[r] & X \ar[rr] \ar@{->>} [rd] && F_0 \ar@{->>}[r] & F\\
&& F_1 \ar@{^{(}->}[ru] && }
\end{equation}
with $F, F_0, F_1 \in \mathcal{F}$.
 Consider the short exact sequences $\mathbb{E}_1: \xymatrix{  F_1 \ar@{^{(}->}[r]   & F_0 \ar@{->>}[r] & F } $ and
 $\mathbb{E}_2 : \xymatrix{ { \rm t}(D) \ar@{^{(}->}[r] & X  \ar@{->>} [r] & F_1 } $ in $\mathcal{C}$.  Since $\Theta_{F_1D}^1$ is an isomorphism,  there exists a short exact sequence $\mathbb{E}'_2$  in $\Ext^1({ \rm q}(F_1), D)$ such that $ \Theta_{F_1D}^1( \mathbb{E}_2')= { \rm t} (\mathbb{E}'_2 ) \eta_{F_1} \equiv \mathbb{E}_2$. So  we have
 $$ \mathbb{E} \equiv \mathbb{E}_2 \circ \mathbb{E}_1 \equiv  ( { \rm t} (\mathbb{E}'_2 ) \eta_{F_1}) \circ \mathbb{E}_1 \equiv  { \rm t} (\mathbb{E}'_2) \circ ( \eta_{F_1}  \mathbb{E}_1).$$
Note that  $ \eta_{F_1}  \mathbb{E}_1$ is a short exact sequence in
$\Ext^1(F, { \rm tq}(F_1))$.
 Similarly, as $\Theta_{F { \rm q }(F_1) }^1$ is an isomorphism,   there exists a short exact sequence $ \mathbb{E}'_1$ in $\Ext^1({ \rm q}(F), { \rm q}(F_1))$  such that ${ \rm t} (\mathbb{E}'_1) \eta_F \equiv \eta_{F_1} \mathbb{E}_1$. Let $\mathbb{E}':= \mathbb{E}'_2 \circ \mathbb{E}'_1$. Notice that $ \mathbb{E}'$ is an exact sequence in $\Ext^2({\rm q}(F), D)$, and
  $${ \rm t}(\mathbb{E}') \eta_F = { \rm t}(\mathbb{E}_2') \circ ({ \rm t}(\mathbb{E}_1') \eta_F) \equiv { \rm t}(\mathbb{E}_2') \circ (\eta_{F_1} \mathbb{E}_1) \equiv \mathbb{E}.$$
Now, we prove that $\Theta_{F D }^2$ is injective.  Let  $\mathbb{E}''$ be an  exact sequence  in $\Ext^2( {\rm q}(F), D)$ such that  $\Theta_{F D  }^2(\mathbb{E}'') ={ \rm t}( \mathbb{E}'') \eta_F$ is trivial. It can be easily verified that there  exist an exact sequence $\mathbb{E}=\mathbb{E}_2 \circ \mathbb{E}_1$  of the form given in  \eqref{exactF} and
a morphism $\mathbb{E} \rightarrow { \rm t}( \mathbb{E}'') $ with  morphisms $\mbox{id}_{ { \rm t} (D) }$ and $\eta_F$ on the most left and right ending objects, respectively, such that  $\mathbb{E}_2$ is  a trivial short exact sequence. Notice that
${ \rm q}(\mathbb{E})=  { \rm q}(\mathbb{E}_2) \circ  { \rm q}(\mathbb{E}_1)$ is exact, and $ {\rm q}(\mathbb{E}_2) $ is a trivial short exact sequence. we consider the composition of morphism of sequences
 $$ { \rm q} (\mathbb{E}) \longrightarrow   { \rm qt}(\mathbb{E}'') \longrightarrow \mathbb{E}'',$$
 where the last morphism is
the canonical counit morphism of the adjunction. The composite has the morphisms $\xi_D$ and $\mbox{id}_{ {\rm q} (F)}$ on the most left and right objects, respectively. So it can be factorized as
$$ { \rm q} (\mathbb{E}) \longrightarrow   \xi_D { \rm q} (\mathbb{E}) \longrightarrow \mathbb{E}''.$$
 $\xi_D { \rm q} (\mathbb{E})$ is a trivial exact sequence which is equivalent to $\mathbb{E}''$. As a consequence, $\Theta_{F }^2$ is an isomorphism for every $F$ in $\mathcal{F}$. By mathematical induction, $\Theta_{F}^n$ is an isomorphism for every $n \geq 1$.

\end{proof}

\subsection{Example.}\label{ex:chain_stalk}
We assume that $\mathbf{C}(\mathcal{C})$ is equipped with the exact structure $\bfC(\mathcal{E})$ (see Example~\ref{ex:chain_complex1}). If $\mathcal{C}$ has (co)kernels, then  for every $i\in \mathbb{Z}$, there is an adjoint triple $( \operatorname{c}_i, \operatorname{s}^i, \operatorname{k}_{i})$
	$$
\xymatrix{ \mathbf{C}(\mathcal{C}) \ar@/^2em/[rrr]^{\mathrm{c}_{i}} \ar@/_2em/[rrr]_{\mathrm{k}_{i}} &&& \mathcal{C} \ar[lll]^{\mathrm{s}^i} },
$$
where $\mathrm{s}^i(C) $ denotes the $i$th sphere chain complex of $C$,  and $\mathrm{c}_i(X)$ and $\mathrm{k}_i(X) $ are the  cokernel and kernel  of differentials $X_{i+1} \longrightarrow X_i$ and $X_i \longrightarrow X_{i-1}$, respectively.  Note that  for every $i\in \mathbb{Z}$, the stalk functor $\mathrm{s}^i$ is exact. So for any chain complex $X$ in $\mathbf{C}(\mathcal{C})$ and $n \geq 1$, there are natural transformations
$$\Theta_X^{n}:\ \Ext^{n}({ \rm c}_i(X), -)  \Longrightarrow \Ext^{n}(X, { \rm s}_i(-))$$
$$\Omega_X^{n}:\ \Ext^{n}(-,{ \rm k}_i(X))  \Longrightarrow \Ext^{n}({ \rm s}_i(-),X)$$
In general, the functors $\mathrm{c}_{i}$ and $\mathrm{k}_{i}$ are not exact (see Example~\ref{ex:chain_stalk2}), and therefore,   $\Theta_X^{n}$   and  $\Omega_X^{n}$ may not be an isomorphism for an arbitrary $X$. Recall that  a chain complex $X$ in $\bfC(\mathcal{C})$ is said to be \textit{acyclic} if  every differential factors as $X_n \twoheadrightarrow \mathrm{k}_{i-1}( X) \hookrightarrow X_{i-1}$  in such a way that each sequence
$$\xymatrix{ \mathrm{k}_i(X) \ar[r] & X_i \ar[r] & \mathrm{k}_{i-1}( X)} $$
is a short exact sequence in $\mcC$; see \cite[Definition 10.1]{Buh}. We let $\mbox{Acyc}$ denote the class of acyclic chain complexes in $\bfC(\mathcal{C})$. By \cite[Proposition 7.5]{Sto13}, we already know  that $\mbox{Acyc}$ satisfies all conditions and their dual forms given in \eqref{subsec:left-right-derived}. As it will be shown in Example~\ref{ex:chain_stalk2}, the functors ${ \rm c}_i$ and ${ \rm k}_i$ restricted on $\mbox{Acyc}$ are exact, and by Lemma~\ref{prop:heredit1},  the natural transformations $\Theta_X^{n}$ and $\Omega_X^{n}$ are isomorphisms for every acyclic chain complex $X$.

\subsection{Lemma.}\label{lemma:derived1}
 Let $\mcF$ and $\mcG$ be classes of objects in $\mathcal{C}$ and $\mathcal{D}$, respectively. For any integer $n \geq 1$, $C \in \mcC$ and $D \in \mathcal{D}$, the followings hold:
\begin{enumerate}[(i)]
	\item Suppose that   ${ \rm t}$ is exact,  $\mathcal{D}$ is an abelian category and $\mcF$ and ${ \rm q}$ satisfy conditions given in \eqref{subsec:left-right-derived}. If for every $1 \leq k \leq n$,
	${ \rm L_k^{ \mcF } q }(C)=0$, then  the natural transformation
$$\Theta_C^{n+1}:\ \Ext^{n+1}({ \rm q}(C), -)  \Longrightarrow \Ext^{n+1}(C, { \rm t}(-))$$
is a monomorphism.

\item Suppose that   ${ \rm q}$ is exact,  $\mathcal{C}$ is an abelian category and $\mcG$ and ${ \rm t}$ satisfies dual of conditions given in \eqref{subsec:left-right-derived}. If for every $1 \leq k \leq n$,
	${ \rm R^k_{ \mcG} t }(D)=0$, then  the natural transformation
$$\Omega_D^{n+1}:\ \Ext^{n+1}(-,{ \rm t}(D))  \Longrightarrow \Ext^{n+1}({ \rm q}(-),D)$$
is a monomorphism.

\end{enumerate}

\begin{proof}
As pointed out in \eqref{subsec:left-right-derived}, using $\mathcal{F}$-resolutions, we have the $k$th derived functor ${ \rm L_k^{\mcF} q }:\ \mathcal{C} \longrightarrow \mathcal{D}$ satisfying ${ \rm L_k^{\mcF}  q }(F)=0$ for every $F \in \mcF$ and $k \geq 1$.

Consider a short exact sequence in $\mathcal{C}$ of the form
\begin{equation}\label{seq1}
\mathbb{E}:\ \xymatrix{ K \ar@{^{(}->}[r] & F \ar@{->>}[r] & C},
\end{equation}
with $F$ in $\mathcal{F}$.    It induces the following long exact sequence  in $\mathcal{D}$
\begin{equation}\label{s2}
\xymatrix{ 0 \ar[r] & { \rm L_1^{\mcF} q}(C) \ar[r] & { \rm q}(K) \ar[r] & { \rm q}(F) \ar[r] & { \rm q}(C) \ar[r] & 0}.
\end{equation}

Suppose that  ${ \rm L_k^{\mcF} q}(C)=0$ for every $1 \leq k \leq n$.  Then   the exact sequence \eqref{s2} turns out to be of  the form
$$\xymatrix{  0 \ar[r] & { \rm q}(K) \ar[r] & { \rm q}(F) \ar[r] &  { \rm q}(C) \ar[r] & 0 }.$$
For a given object  $D$  in $\mathcal{D}$, applying  the contravariant functor $\Hom(-,D)$, and subsequently, using the adjoint pair  $({ \rm q},{ \rm t})$ and Lemma \ref{prop:heredit1}, we have the following commutative diagram of exact sequences
$$
\scalebox{0.9}{\xymatrix{
 \Ext^{n}({ \rm q}(F), D) \ar[r] \ar[d]^{\Theta^{n}_{FD}}_{\cong} &  \Ext^{n}({ \rm q}(K), D) \ar[r] \ar[d]^{\Theta^{n}_{KD}} & \Ext^{n+1}({ \rm q}(C), D) \ar[r] \ar[d]^{\Theta^{n+1}_{CD}} & \Ext^{n+1}({ \rm q}(F), D)  \ar[d]^{\Theta^{n+1}_{FD}}_{\cong} \\
 \Ext^{n}(F, { \rm t}(D)) \ar[r]  & \Ext^{n}(F, { \rm t}(D)) \ar[r]  & \Ext^{n+1}(C, { \rm t}(D)) \ar[r]  & \Ext^{n+1}(F, { \rm t}(D)) }}
$$

Now, we apply mathematical induction on $n$. If $n=1$, then by Lemma \ref{lemma:ext_adjoint_pair},  $\Theta^{1}_{KD}$ is a monomorphism, and therefore, so is  $\Theta^{2}_{CD}$.

Now, suppose that the statement is true for some $n \geq 1$, and ${ \rm L_k^{\mcF} q}(C)=0$ for every $1 \leq k \leq n+1$. Then ${ \rm L_k^{\mcF} q}(K)=0$ for every $1 \leq k \leq n$. By induction hypothesis,  $\Theta^{n+1}_{KD}$ is a monomorphism, then so is $\Theta^{n+2}_{CD}$

\end{proof}

\subsection{Proposition.}\label{lemma:derived_zero}
 Let $\mcF$ and $\mcG$ be classes of objects in $\mathcal{C}$ and $\mathcal{D}$, respectively. For any integer $n \geq 1$, $C \in \mcC$ and $D \in \mathcal{D}$, the followings hold:
\begin{enumerate}[(i)]
	\item If  ${ \rm t}$ is exact,  $\mathcal{D}$ is an abelian category, and $\mcF$  and ${ \rm q}$ satisfy conditions given in \eqref{subsec:left-right-derived}, then the following are equivalent:
	\begin{enumerate}[(a)]
	\item  For every $1 \leq k \leq n$,   ${ \rm L^{\mcF} _k q }(C)=0$;
\item For every $1 \leq k \leq n$, the natural transformation
$$\Theta_C^k:\ \Ext^k({ \rm q}(C), -)  \Longrightarrow \Ext^k(C, { \rm t}(-)),$$
is an isomorphism.

\end{enumerate}
\item  If ${ \rm q}$ is exact,  $\mathcal{C}$ is an abelian category, and $\mathcal{G}$ and ${ \rm t}$ satisfy dual conditions given in \eqref{subsec:left-right-derived}, then the  followings hold:
\begin{enumerate}[(a)]
\item For every $1 \leq k \leq n$,  ${ \rm R^k_{\mcG} t}(D)=0$
\item  For every $1 \leq k \leq n$,  the natural transformation
$$\Omega_D^k:\ \Ext^k(-, { \rm t}(D))  \Longrightarrow \Ext^k({ \rm q}(-), D),$$
is an isomorphism.
\end{enumerate}
\end{enumerate}

\begin{proof}
Consider a short exact sequence in $\mathcal{C}$ of the form
\begin{equation}\label{s1}
\mathbb{E}:\ \xymatrix{ K \ar@{^{(}->}[r] & F \ar@{->>}[r] & C},
\end{equation}
with $F$ in $\mathcal{F}$.   It induces the following long exact sequence  in $\mathcal{D}$
\begin{equation}\label{seq2}
\xymatrix{ 0 \ar[r] &{ \rm L_1^{\mcF} q}(C) \ar[r] & { \rm q}(F) \ar[r] & { \rm q}(P) \ar[r] & { \rm q}(C) \ar[r] & 0}.
\end{equation}

 We apply mathematical induction on $n$ to prove the equivalence  of the statements (a) and (b). If  ${ \rm L_1^{\mcF} q}(C)=0$, then    the exact sequence \eqref{seq2} is turned out being of   the form
$$\xymatrix{  0 \ar[r] & { \rm q}(K) \ar[r] & { \rm q}(F) \ar[r] &  { \rm q}(C) \ar[r] & 0 },$$
which by  Lemma \ref{prop:heredit1}, yields
 the following commutative diagram of exact sequences of functors
$$
\scalebox{0.8}{\xymatrix{
 \Hom({ \rm q}(F), -) \ar@2[r] \ar@2[d]_{\cong} &  \Hom({ \rm q}(K), -) \ar@2[r] \ar@2[d]_{\cong} & \Ext^1({ \rm q}(C), -) \ar@2[r] \ar@2[d]^{\Theta^1_{C}} & \Ext^1({ \rm q}(F), -) \ar@2[r] \ar@2[d]^{\Theta^1_{F}}_{\cong} & \Ext^1({ \rm q}(K), -) \ar@2[d]^{\Theta^1_{K}}\\
 \Hom(F, { \rm t}(-)) \ar@2[r]  &  \Hom(F, { \rm t}(-)) \ar@2[r] ] & \Ext^1(C, { \rm t}(-)) \ar@2[r]  & \Ext^1(F, { \rm t}(-)) \ar@2[r]  & \Ext^1(K, { \rm t}(-))
}}
$$
Using the fact that $\Theta^1_{CD}$ and $\Theta^1_{KD}$ are monomorphisms, one can easily show that $\Theta^1_{CD}$ is in fact an isomorphism.

\sloppy For the converse, suppose that the natural transformation
$${\Theta^1_C: \Ext^1({ \rm q}(C), -) \Longrightarrow \Ext^1(C, { \rm t}(-))}$$ is an isomorphism. By Proposition \ref{prop:adjoints}, ${ \rm q}(\mathbb{E})$ is a short exact sequence in $\mathcal{D}$, and hence, ${ \rm L_1^{\mcF} q}(C)=0$.


Now, suppose that the equivalence (a)$\Leftrightarrow$(b) is true for any $n \geq 1$ and any object in $ \mathcal{C}$.  Assume that   ${ \rm L_k^{\mcF} q}(C)=0$ for every $1 \leq k \leq n+1$. By  induction step,  the natural transformation
$$\Theta^k_C:\ \Ext^k({ \rm q}(C), -) \Longrightarrow \Ext^k(C, { \rm t}(-))$$
is an isomorphism for every $1 \leq k \leq n$, and by Lemma  \ref{lemma:derived1}, $\Theta^{n+1}_C$ is a monomorphism. And from the long exact sequence of derived functors ${ \rm L^{\mcF} _kq}(-)$ obtained through the short exact sequence $\mathbb{E}$, we have
$$0={ \rm L_{k+1}^{\mcF}  q}(C) \cong { \rm L_{k} ^{\mcF} q}(K)$$
for every $1 \leq k \leq n$. By induction step, $\Theta_K^k$ is an isomorphism for every $1 \leq k \leq n$, and  by Lemma  \ref{lemma:derived1}, $\Theta_K^{n+1}$ is a monomorphism.

As   ${ \rm L_1^{\mcF} q}(C)=0$,  ${ \rm q}(\mathbb{E})$ is a   short exact sequence in $\mathcal{D}$, which induces the following commutative diagram
$$
\scalebox{0.9}{\xymatrix{
  \Ext^{n}({ \rm q}(K), -) \ar@2[r] \ar@2[d]^{\Theta^{n}_{K}}_{\cong}& \Ext^{n+1}({ \rm q}(C), -) \ar@2[r] \ar@2[d]^{\Theta^{n+1}_{C}} & \Ext^{n+1}({ \rm q}(F), -)  \ar@2[d]^{\Theta^{n+1}_{F}}_{\cong}  \ar@2[r] &  \Ext^{n+1}(K, { \rm t}(-)) \ar@2[d]^{\Theta^{n+1}_{K}}\\
  \Ext^{n}(F, { \rm t}(-)) \ar@2[r] ] & \Ext^{n+1}(C, { \rm t}(-)) \ar@2[r]  & \Ext^{n+1}(F, { \rm t}(-))  \ar@2[r] &  \Ext^{n+1}(K, { \rm t}(-)) }}
$$
with exact rows. From the diagram, it is easy to verify that $\Theta^{n+1}_C$ is an epimorphism.

For the converse,  suppose that the natural transformation $\Theta_C^k$ is an isomorphism  for every $1 \leq k \leq n+1$. By induction hypothesis, ${ \rm L_k^{\mcF} q}(C)=0$  for every $1 \leq k \leq n$.
As a particular case, ${ \rm L_1^{\mcF} q}(C)=0$, and therefore, ${ \rm q}(\mathbb{E})$ is  a short exact sequence in $\mathcal{D}$. Applying five lemma, one can show that $\Theta_K^k $ is an isomorphism for every $1 \leq k \leq n$. Again, by induction hypothesis,  ${ \rm L_k^{\mcF} q}(K)=0$ for every $1 \leq k \leq n$. By shifting dimension, we have
${ \rm L_{n+1}^{\mcF} }{ \rm q}(C) \cong { L_{n}^{\mcF} }{ \rm q}(K)=0$.

\end{proof}

\subsection{Example.}\label{ex:adj_ring_ext1}
Let $f:\ R \rightarrow S$  a  ring homomorphism. Consider the  adjoint pair
$$
\xymatrix{   R \Mod  \ar@/^1em/[rrr]^{ S \otimes_R -}   &&&  S \Mod  \ar@/^1em/[lll]-^{\operatorname{res}}},
$$
where $\mbox{res}:\ S \Mod \longrightarrow R \Mod$ denotes the functor of restriction of scalars through  the ring homomorphism $f$. Note that $\operatorname{res}$ is an exact functor, and
$${{ \rm L_k}( S \otimes_R-)=  \mbox{Tor}_R^k(S,-)}.$$
By Proposition   \ref{lemma:derived_zero}-(i), for a left $R$-module $M$, $\operatorname{Tor}_R^k(S,M)=0$ for every $1 \leq k \leq n$ if and only if $\Theta_M^k:\ \Ext^k( S \otimes_R M, -)  \Longrightarrow \Ext^k(M, \operatorname{res}(-)),$
is an isomorphism for every $1 \leq k \leq n$.

Similarly, for the adjoint pair
$$
\xymatrix{   S\Mod  \ar@/^1em/[rrr]^{ \operatorname{res}}   &&& R \Mod  \ar@/^1em/[lll]^{ \Hom(S,-)}      }
$$
By Proposition \ref{lemma:derived_zero}, for a left $R$-module $N$, $\Ext^k(S,N)=0$ for every $1 \leq k \leq n$ if and only if $\Omega_N^k:\ \Ext^k( -, \Hom(S,N))  \Longrightarrow \Ext^k( \operatorname{res}(-), N)$ is an isomorphism for every $1 \leq k \leq n$.

\section{Cotorsion pairs induced from adjoint pairs}
In this section, we apply the results from   Section 2 to  investigate cotorsion pairs arising through  an adjoint pair. We also provide  characterizations of  objects in the right/left perpendicular classes that cover  well-known examples.

We keep assuming that $\mathcal{C}$ and $\mathcal{D}$ are exact categories. For a given family   $\{ x_i\}_{i\in I}$  of  functors from $\mathcal{C}$ to $\mathcal{D}$, and classes  $\mcF$ and $\mcG$ of objects in $\mathcal{C}$ and $\mathcal{D}$, respectively,  we let
\begin{align*}
x_*(\mcF)&:=\{ x_i(F)  \mid\ F \in \mcF \textrm{ and } i \in I \}\\
x_*^{-1}(\mcG)&:=\{ C \in \mathcal{C} \mid\ x_i(C) \in \mcG, \textrm{ for every } i \in I \}.
\end{align*}
An object in $ \mathcal{C}$ will be said   to be $x_*$-(co)flat if it is $x_i$-(co)flat for every $i \in I$.
\subsection{Theorem.}\label{Prop:degrewise-cot.pair}
Let $\{({ \rm q}_i, { \rm t}_i)\}_{i \in I}$ be a family of adjoint pairs
	$$
\xymatrix{\mathcal{C} \ar@/^1em/[rrr]^{{ \rm q}_i} &&& \mathcal{D} \ar@/^1em/[lll]^{{ \rm t}_i} }.
$$
For given cotorsion pairs $(\mcF, \mcG)$ and $(\mcF',\mcG')$  in $\mcC$ and  $\mathcal{D}$, respectively, the following statements hold.
\begin{enumerate}[(i)]
	\item Suppose that for every $i \in I$, every object in $\mathcal{F}$ is ${ \rm q}_i$-flat.
	\begin{enumerate}[(a)]
	\item If for every $i \in I$,  ${ \rm t}_i$  preserves short exact sequences in $\Ext^1({ \rm q}_i( \mcF), -)$, then
	 $$({ \rm q}_*(\mcF))^\perp= { \rm t}_*^{-1} (\mcG).$$
	
	 \item If $\mcF$ is a generating class in $\mcC$, then
	 $$({ \rm q}_*(\mcF))^\perp=\{ D \in { \rm t}_*^{-1} (\mcG) \mid\   D \textrm{ is a }{ \rm t}_*\textrm{-coflat object}  \}. $$
	
	 \end{enumerate}
\item Suppose that for every $i \in I$, every objeect in $\mcG'$ is ${ \rm t}_i$-coflat. 

 \begin{enumerate}[(a)]
 \item If for every $i \in I$,  ${ \rm q}_i$ preserves short exact sequences in $\Ext^1(-, { \rm t}_i( \mcG'))$, then
$$ {}^\perp ({ \rm t}_*(\mcG')) =  { \rm q}_*^{-1}(\mcF').$$

\item If $\mcG'$ is a cogenerating class in $\mathcal{D}$, then
$${}^\perp ({ \rm t}_*(\mcG'))= \{ C \in { \rm q}_*^{-1}(\mcF')  \mid\ C \textrm{ is a }{ \rm q}_*\textrm{-flat object}  \}.$$
\end{enumerate}
\end{enumerate}

\begin{proof}
We only prove the first statement, as the other can be proved by dual arguments. Firstly, as $\Omega^1_{F}$ is a monomorphism for every $F \in \mathcal{F}$ by Lemma \ref{lemma:ext_adjoint_pair}-(ii),
$$({ \rm q}_*(\mcF))^\perp \subseteq { \rm t}_*^{-1} (\mcG).$$

If for every $i \in I$,  ${ \rm t}_i$  preserves short exact sequences in $\Ext^1({ \rm q}_i( \mcF), -)$, then by Lemma \ref{lemma:ext_adjoint_pair}-(iii),
	 $$({ \rm q}_*(\mcF))^\perp= { \rm t}_*^{-1} (\mcG).$$

On the other hand, if $\mcF$ is a generating class, then by Lemma \ref{lemma:exactt}-(i), any object in $({ \rm q}_*(\mcF))^\perp$ is a ${ \rm t}_i$-coflat object for every $i \in I$. The converse inclusion $\supseteq$ in  (i)-(b) is immediate   by Lemma \ref{lemma:ext_adjoint_pair}-(i).

\end{proof}

\subsection{Remark.}\label{remark:clases-HJ} Note that if for every $i \in I$, ${ \rm q}_i$ is exact, $\mathcal{C}$ is an abelian category, and the $k$th right derived functor ${ \rm R_{\mcG}^kt_i}$ exists as described in \eqref{subsec:left-right-derived}, then by Proposition \ref{prop:adjoints}-(ii) and Proposition~\ref{lemma:derived_zero}-(ii), for every $i \in I$,  an object $D$ in $\mathcal{D}$ is ${ \rm t}_i$-coflat   is equivalent to ${ \rm R^1_{\mcG}t}_i(D)=0$.  In particular, if $\mathcal{D}$ is an abelian category with enough injectives, and $\mcC$ has enough projectives then   the classes $\Psi(\mathcal{G})$ and $\Phi(\mathcal{F'})$ defined in \cite[Definition 1.5]{HJ19} are exactly the classes  ${ \rm q}_*(\mathcal{F})^\perp$ and  ${}^\perp { \rm t}_*(\mathcal{G'})$ described  in Theorem \ref{Prop:degrewise-cot.pair}-(i)-(b) and (ii)-(b), respectively.

\subsection{Corollary.}\label{setup:3}
 Let $\{({ \rm q}_i, { \rm t}_i, { \rm p_i})\}_{i \in I}$ be a family of adjoint triples
	$$
\xymatrix{\mathcal{C} \ar@/^2em/[rrr]^{{ \rm q}_i} \ar@/_2em/[rrr]_{ { \rm p_i}} &&& \mathcal{D} \ar[lll]^{{ \rm t}_i} }.
$$
For  given cotorsion pairs $(\mcF,\mcG)$ and $(\mcF',\mcG')$ in $\mathcal{C}$ and $\mathcal{D}$, respectively, the following hold:
\begin{enumerate}[(i)]
\item  Suppose that for every $i\in I$, ${ \rm q}_i$ preserves short exact sequences in $\Ext^1( \mcF,{ \rm t}_i(-))$.
\begin{enumerate}[(a)]
\item  If for every $i \in I$, ${ \rm t}_i$ preserves short exact sequences in $\Ext^1({ \rm q}_i( \mcF),-)$, then
$$ { \rm q}_*(\mathcal{F})^\perp= { \rm t}_*^{-1} (\mcG);$$
\item If  $\mcF$ is a generating class in $\mathcal{C}$, then
$$ { \rm q}_*(\mathcal{F})^\perp=\{ D \in { \rm t}_*^{-1} (\mcG) \mid\  D \textrm{ is a }{ \rm t}_*\textrm{-coflat object} \}.$$
\end{enumerate}
\item Suppose that  for every $i\in I$,  ${ \rm p_i}$ preserves short exact sequences in $\Ext^1({ \rm t}_i(-), \mathcal{G})$.
\begin{enumerate}[(a)]
\item If for every $i \in I$, ${ \rm t}_i$ preserves short exact sequences in $\Ext^1(-, { \rm p_i}( \mathcal{G}))$, then
$$ {}^\perp { \rm p}_*(\mathcal{G} )=  { \rm t}_*^{-1}(\mathcal{F});$$
\item If $\mcG$ is a cogenerating class in $\mathcal{C}$, then
$$ {}^\perp { \rm p}_*(\mathcal{G} )= \{ D \in { \rm t}_*^{-1}(\mcF)  \mid\  D \textrm{ is a }
{ \rm t}_*\textrm{-flat object} \}. $$
\end{enumerate}

\item Suppose that  for every $i\in I$, ${ \rm t}_i$  preserves short exact sequences in $\Ext^1( \mathcal{F}',{ \rm p}_i(-))$.
\begin{enumerate}[(a)]
 \item If for every $i\in I$, ${ \rm p}_i$ preserves short exact sequences in $\Ext^1( { \rm t}_i (\mathcal{F}'),-)$, then
$$ { \rm t}_*(\mathcal{F}')^\perp= { \rm p}_*^{-1}(\mcG').$$

\item If $\mathcal{F}'$ is a generating class in $\mathcal{D}$, then
$$ { \rm t}_*(\mathcal{F}') ^\perp=\{   C \in { \rm p}_*^{-1}(\mcG')  \mid\  C \textrm{ is a }{ \rm p}_* \textrm{-coflat object} \}.$$

\end{enumerate}

\item
 Suppose that for every $i \in I$, ${ \rm t}_i$ preserves short exact sequences in  $\Ext^1({ \rm q}_i(-),\mcG')$.
\begin{enumerate}[(a)]
 \item If for every $i \in I$,  ${ \rm q}_i$ preserves short exact sequences in $\Ext^1(-, { \rm t}_i( \mcG'))$, then
$$ {}^\perp { \rm t}_*(\mcG') =  { \rm q}_*^{-1}(\mcF').$$

\item If $\mcG'$ is a cogenerating class in $\mathcal{D}$, then
$${}^\perp { \rm t}_*(\mcG')= \{ C \in { \rm q}_*^{-1}(\mcF')  \mid\ C \textrm{ is a }{ \rm q}_*\textrm{-flat object}  \}.$$
\end{enumerate}

\end{enumerate}
\begin{proof}
Apply Theorem \ref{Prop:degrewise-cot.pair}.
\end{proof}

\subsection{Remark.}Since for every $i \in I$, the middle functor ${ \rm t}_i$ given in Corollary \ref{setup:3} is both right and left adjoint functor, it preserves all (co)limits, but it is not necessarily an exact functor between exact categories. However, it is exact whenever $\mathcal{C}$ is an abelian category.
\subsection{Example.}\label{ex:chain_ev_cot}
Consider  Example  \ref{ex:chain_complex1}.  For a given cotorsion pair $(\mcF, \mcG)$ in $\mcC$ and $i \in \Z$, we let
\begin{align*}
\mathbf{C}(\mcF)_i&:=\{ X \in  \mathbf{C}(\mcC) \mid\ X_i  \in \mcF  \};\\
\mathbf{C}(\mcG)_i&:=\{ X \in  \mathbf{C}(\mcC) \mid\ X_i  \in \mcG  \}.
\end{align*}
As $\mbox{D}^i$ and $\ev_i$ are exact functors, by Corollary~\ref{setup:3}, there are two  cotorsion pairs in $\mathbf{C}(\mcC)$
$$(  {}^\perp \mathbf{C}(\mcG)_i, \mathbf{C}(\mcG)_i) \quad \textrm{ and } \quad (\mathbf{C}(\mcF)_{i}, \mathbf{C}(\mcF)_{i}^\perp )$$
generated by $\mbox{D}^i(\mathcal{F})$ and cogenerated by $\mbox{D}^{i+1}(\mathcal{G})$, respectively. If we consider all $i \in \Z$, we have the following cotorsion pairs
$$(  {}^\perp \mathbf{C}(\mcG), \mathbf{C}(\mcG)) \quad \textrm{ and } \quad (\mathbf{C}(\mcF), \mathbf{C}(\mcF)^\perp );$$
see \cite[Proposition 3.2]{Gil08}.

If $(\mcF', \mcG')$ is a cotorsion pair in $\mathbf{C}(\mcC)$, then
\begin{align*}
\mbox{ev}_*( \mathcal{F}')^\perp & = \{ C \in \mathcal{C} \mid\ \mbox{D}^i(C) \in \mathcal{G}' \textrm{ for every } i \in \Z \};\\
^\perp \mbox{ev}_*( \mathcal{G}') & = \{ C \in \mathcal{C} \mid\ \mbox{D}^i(C) \in \mathcal{F}' \textrm{ for every } i \in \Z \};
\end{align*}
see \cite[Lemma 3.5]{Gil08}.

\subsection{Example.}\label{ex:chain_stalk2}
Consider  Example~\ref{ex:chain_stalk}. Let $(\mathcal{F}, \mathcal{G})$ be a cotorsion pair in $\mcC$. Note that $\mathrm{k}_i$ doesn't preserves all short exact sequences in $\Ext^1( \mathrm{s}^i(\mathcal{F}),-)$ for any $i \in \Z$; for example, consider the canonical short exact sequence
$$ \xymatrix{ \mathrm{s}^{i-1}(F) \ar@{^{(}->}[r] & \operatorname{D}^i(F) \ar@{->>}[r] &  \mathrm{s}^i(F) },$$
where $F$ is a non-zero object in $\mathcal{F}$. On the other hand, if $\mathcal{F}$ is a generating class in $\mathcal{C}$, then
$$\mathrm{s}^i(\mathcal{F})^\perp=\{X \in \mathbf{C}(\mathcal{C}) \mid\  \mathrm{k}_i(X) \in \mathcal{G} \textrm{ and } X \textrm{ is a  } \mathrm{k}_i\textrm{-coflat object}\}.$$
We now show that a chain complex  $X$ in $\mathbf{C}(\mathcal{C}) $ is a $\mathrm{k}_i$-coflat object  if and only if the canonical sequence
$$\xymatrix{ \mathrm{k}_i(X) \ar[r] & X_i \ar[r] & \mathrm{k}_{i-1}( X)} $$
is a short exact sequence in $\mcC$.

Suppose that  $X$  is a $\mathrm{k}_i$-coflat object for some $i \in \Z$, then consider the canonical morphism $f_{i-1}: \mathrm{s}^{i-1}( \mathrm{k}_{i-1}(X)) \longrightarrow X$, which induces the degreewise split short exact sequence in $\mathbf{C}(\mathcal{C})$
$$\mathbb{E}:\ \xymatrix{X \ar@{^{(}->}[r] & \operatorname{cone}(f_{i-1}) \ar@{->>}[r] & \mathrm{s}^{i}( \mathrm{k}_{i-1}(X))},$$
where $\operatorname{cone}(f_{i-1})$ is the cone of the morphism $f_i$. By assumption, $\mathrm{k}_{i} (\mathbb{E})$ is a short exact sequence in $\mathcal{C}$, and it is of the form
$$\mathrm{k}_{i} (\mathbb{E}):\ \xymatrix{\mathrm{k}_{i}(X) \ar@{^{(}->}[r] & X_{i} \ar@{->>}[r] & \mathrm{k}_{i-1}(X)}.$$

Conversely, suppose that the canonical sequence
$$\mathrm{k}_{i} (\mathbb{E}):\ \xymatrix{\mathrm{k}_{i}(X) \ar@{^{(}->}[r] & X_{i} \ar@{->>}[r] & \mathrm{k}_{i-1}(X)}$$
 is a short exact sequence in $\mcC$.
 By Lemma \ref{lemma:lem}, it is enough to show that $\mathrm{k}_i$ preserves all short exact sequences in $\Ext^1(\mathrm{s}^i(-),X)$. Consider a short exact sequence in $\mathbf{C}(\mathcal{C})$ of the form
$$\mathbb{E}:\ \xymatrix{X \ar@{^{(}->}[r] & Y \ar@{->>}[r] & \mathrm{s}^{i}(C)}.$$
Without lost of generality, we can assume that $X_j=Y_j$ for every $i \neq j \in \Z$.
The short exact sequence $\mathbb{E}$ induces the following commutative diagram
$$
\xymatrix{
X_i \ar@{^{(}->}[rrr] \ar[dd] \ar@{->>}[rd] &&& Y_i  \ar[ld] \ar@{->>}[r] \ar[dd] & C \ar[dd] \\
&\mathrm{k}_{i-1}(X) \ar[ld] \ar@{=}[r]&\mathrm{k}_{i-1}(Y) \ar[rd] &&\\
X_{i-1} \ar@{=}[rrr] \ar[d] &&& X_{i-1}  \ar[r]\ar[d]& 0\ar@{=}[d]\\
X_{i-2} \ar@{=}[rrr] &&& X_{i-2} \ar[r]& 0
}
$$
By assumption, $\mathcal{C}$ has kernels, and therefore, the morphism $Y_i \longrightarrow \mathrm{k}_{i-1}(Y)$ is in fact an admissible epic.
Applying Snake Lemma  \cite[Corollary 8.13]{Buh} to the diagram
$$
\xymatrix{
X_i \ar@{^{(}->}[r] \ar[d] \ar@{->>}[d] & Y_i  \ar@{->>}[d] \ar@{->>}[r]  & C \ar[d] \\
\mathrm{k}_{i-1}(X) \ar@{=}[r]&\mathrm{k}_{i-1}(Y) \ar[r] & 0
}
$$
 we obtain the short exact sequence in $\mathcal{C}$
$$\mathrm{k}_{i} (\mathbb{E}):\ \xymatrix{ \mathrm{k}_{i}(X) \ar@{^{(}->}[r] & \mathrm{k}_{i}(Y) \ar@{->>}[r] &C=\mathrm{k}_{i}( \mathrm{s}^{i}(C))  }.$$
If we consider all integers $i$, we have
\begin{align*}
\mathrm{s}^*(\mathcal{F})^\perp &=\{X \in \mathbf{C}(\mathcal{C}) \mid\ \mathrm{k}_*(X) \in \mathcal{G} \textrm{ and } X \textrm{ is a } \mathrm{k}_* \textrm{-coflat object} \}\\
&= \{ X \in \mathbf{C}(\mathcal{C}) \mid\ \mathrm{k}_*(X) \in \mathcal{G} \textrm{ and } X \textrm{ is acyclic}\}
\end{align*}

Similarly, $\mathrm{c}_i$ doesn't preserves all short exact sequences in $\Ext^1(-, \mathrm{s}^i(\mathcal{G}))$
for any $i \in \Z$, and  if  $\mathcal{G}$ is a cogenerating class in $\mcC$, then
\begin{align*}
^\perp \mathrm{s}^*(\mathcal{G}) &=\{X \in \mathbf{C}(\mathcal{C}) \mid\ \mathrm{c}_*(X) \in \mathcal{F} \textrm{ and } X \textrm{ is a } \mathrm{c}_* \textrm{-flat object} \}\\
&= \{ X \in \mathbf{C}(\mathcal{C}) \mid\ \mathrm{c}_*(X) \in \mathcal{F} \textrm{ and } X \textrm{ is acyclic}\}.
\end{align*}

\subsection{Example.}\label{example:character_dual}
Consider the following adjoint triple
$$
\xymatrix{ (R \Mod)^{\textrm{op}}    \ar@/_2em/[rrr]_{\Hom(-,\mathbb{Q}/\Z)}  \ar@/^2em/[rrr]^{\Hom(-,\mathbb{Q}/\Z)}  &&& \mbox{Mod-}R \ar[lll]_{\Hom(-,\mathbb{Q}/\Z)} },
$$
 Since $\mathbb{Q}/\Z$ is an injective abelian group,     $\Hom(-,\mathbb{Q}/\Z)$ is an exact functor. If $(\mcF, \mcG)$ is  a cotorsion pair in $R \Mod$, then
$$\Hom(\mathcal{G},\mathbb{Q}/\Z)^\perp=\{  M \in  \mbox{Mod-}R \mid\  \Hom(M,\mathbb{Q}/\Z) \in \mathcal{F} \} $$
$${}^\perp \Hom(\mathcal{F},\mathbb{Q}/\Z)=\{  M \in  \mbox{Mod-}R  \mid\  \Hom(M,\mathbb{Q}/\Z) \in \mathcal{G} \}. $$

\subsection{Example.}
Consider the degreewise pure and split exact structures on $\mathbf{C}(R)$ and $\mathbf{C}(R^{\textrm{op}})$, respectively. The adjoint triple of  Example \ref{example:character_dual} can be extended to   the following  adjoint triple of functoors
$$
\xymatrix{ \mathbf{C}(R)^{\textrm{op}}   \ar@/_2em/[rrr]_{\Hom(-,\mathbb{Q}/\Z)}  \ar@/^2em/[rrr]^{\Hom(-,\mathbb{Q}/\Z)}  &&& \mathbf{C}(R^{\textrm{op}}) \ar[lll]_{\Hom(-,\mathbb{Q}/\Z)} },
$$
which are exact.
We let $\mbox{Acyc}(R^{\textrm{op}})$ and $\mbox{K-Inj}(R^{\textrm{op}})$ denote the classes of acyclic and k-injective chain complexes of $R^\textrm{op}$-modules: see \cite{Spa88}. It can be easily verifed that   $( \mbox{Acyc}(R^{\textrm{op}}), \mbox{K-Inj}(R^{\textrm{op}}))$ is a cotorsion pair  in $\mathbf{C}(R^{\textrm{op}})$.  So we have
\begin{align*}
{}^\perp \Hom( \mbox{Acyc}(R^{\textrm{op}})  ,\mathbb{Q}/\Z)&=\{  X \in  \mathbf{C}(R)  \mid\  \Hom(X,\mathbb{Q}/\Z) \in \mbox{K-Inj}(R^{\textrm{op}}) \}\\
&=\{  X \in  \mathbf{C}(R)  \mid\  X \textrm{ is k-flat} \}
\end{align*}
The last equality follows from the isomorphism
$$\operatorname{H}_i(\Hom( Y \otimes_R X,\mathbb{Q}/\Z) )\cong \Hom_{\mathbf{K}(R^{\textrm{op}})}(Y, \Hom(X,\mathbb{Q}/\Z)),$$
where $Y \in \bfC(R^{\textrm{op}} )$ and $Y \otimes_R X$ denotes the total tensor product of chain complexes.

\subsection{Example.}\label{ex:1} Let $\{f_i: R \rightarrow S_i\}_{i \in I}$ be a family of ring homomorphisms. For every $i \in I$, consider the adjoint triple  $( S_i \otimes_R -, \operatorname{res}_i,  \Hom_R(S_i,-) )$ as given in Example~\ref{ex:adj_ring_ext1},
 where $\mbox{res}_i:\ S_i \Mod \longrightarrow R \Mod$ denotes the functor of restriction of scalars through $f_i$. They induce a family of adjoint triples
$$
\xymatrix{   R\Mod  \ar@/^2em/[rrr]^{ \iota_i \circ (S_i \otimes_R -)}   \ar@/_2em/[rrr]_{ \iota_i \circ \Hom_R(S_i,-) }   &&& \prod_{i\in I} S_i \Mod  \ar[lll]_{\operatorname{res}_i \circ \pi_i}      }.
$$
Note that $\operatorname{res}_i \circ \pi_i$ is an exact functor. 
Let  $(\mcF, \mcG)$ be a cotorsion pair  in $R \Mod$.
\begin{enumerate}[(i)]
	
	\item By Lemma \ref{lemma:derived_zero}, for every $i \in I$,
	$ (S_i \otimes_R -)$ preserves short exact sequences in $\Ext^1( \mcF, \operatorname{res}_i(-))$ if and only if  $\mbox{Tor}_R^1(S_i, F)=0$ for every $F \in \mathcal{F}$. In such case,
	$$(\iota_* \circ ( S_* \otimes_R \mathcal{F}))^\perp=\{ (N_i)_{i \in I} \in \prod_{i\in I} S_i \Mod \mid  \operatorname{res}_i(N_i) \in \mathcal{G} \textrm{ for every } i \in I\}.$$

\item By Lemma \ref{lemma:derived_zero}, for every $i \in I$  the functor $\Hom_R(S_i,-)$ preserves short exact sequences in $\Ext^1(\operatorname{res}_i(-), \mcG)$ if and only if $\Ext_R^1(S_i,G)=0$ for every $G \in \mathcal{G}$, that is, $ \operatorname{res}_i(S_i) \in \mathcal{F}$. In such case,
$${}^\perp(\iota_* \circ \Hom_R(S_*,\mathcal{G}))=\{  (N_i)_{i \in I} \in \prod_{i\in I} S_i \Mod \mid  \operatorname{res}_i(N_i) \in \mathcal{F} \textrm{ for every } i \in I \}. $$

\end{enumerate}

On the other hand, by Lemma \ref{lemma:product_cot_pair}, any cotorsion pair in $\prod_{i\in I} S_i \Mod$ is just a family $\{ (\mcF_i, \mcG_i) \}_{i \in I}$ of cotorsion pairs $( \mcF_i, \mcG_i)$ in $S_i \Mod$.
\begin{enumerate}[(i)]
\item As for every $i \in I$, $\mcF_i$ is a generating class in $S_i\Mod$,
$$\mbox{res}_*( \mcF_*)^\perp= \{N \in R \Mod \mid \ \Hom_R(S_i,N) \in \mcG_i, \Ext^1_R(S_i,N)=0 \textrm{ for every } i \in I \}$$
\item As for every $i \in I$, $\mcG_i$ is a cogenerating class in $S_i\Mod$,
$$^\perp \mbox{res}_*( \mcG_*)= \{N \in R \Mod \mid \ S_i \otimes_R N \in \mcF_i, \mbox{Tor}_R^1(S_i,N)=0 \textrm{ for every } i \in I \}$$
\end{enumerate}

\subsection{Example.}
We know consider a very particular case of Example \ref{ex:1}. Let $R$ be a commutative ring. For a given element  $r$ of $R$, consider the localization of $R$ at $r$, denoted by $R[r^{-1}]$, together with the  ring homomorphism $R \rightarrow R[r^{-1}]$. As $R[r^{-1}]$ is a flat $R$-$R$-bimodule,  the functor $R[r^{-1}] \otimes_R - $ is exact, and  therefore, $\operatorname{Tor}_R^1(R[r^{-1}] ,-)=0$.  Note that for any $M \in R\Mod$,
$$R[r^{-1}] \otimes_R N  \cong N[r^{-1}].$$
 For a given cotorsion pair  $(\mcF_r,\mcG_r)$ in $R[r^{-1}] \Mod$, we have
\begin{equation*}
\begin{split}
{}^\perp ( \operatorname{res} (\mcG_r))&=  \{ N \in R \Mod \mid N[r^{-1}] \in  \mcF_r  \}; \\
( \operatorname{res}(\mathcal{F}_r) )^\perp &=\{   N \in R \Mod  \mid  \Hom_R(R[r^{-1}], N) \in \mcG_r \textrm{ and } \Ext^1_R(R[r^{-1}],N)=0 \}.
\end{split}
\end{equation*}
As a particular case, consider the trivial cotorsion pair $(\Proj_r, R[r^{-1}] \Mod)$ in $R[r^{-1}] \Mod$, where $\Proj_r$ denotes the class of projective $R[r^{-1}]$-modules. So we have
\begin{equation*}
\begin{split}
{}^\perp (\operatorname{res}(R[r^{-1}] \Mod ))&=\{ N \in R\Mod \mid   N[r^{-1}]\in \Proj   \}\\
(\operatorname{res}( \Proj_r))^\perp
&=\{ N \in R\Mod \mid \mbox{Ext}_R^1(R[r^{-1}] ,N)=0   \}\\
&=\{ r\textrm{-contraadjusted $R$-modules} \}
\end{split}
\end{equation*}
Consider the family of adjoint triples
$$
\xymatrix{   R\Mod  \ar@/^2em/[rrr]^{ \iota_r \circ (R[r^{-1}] \otimes_R -)}   \ar@/_2em/[rrr]_{ \iota_r \circ \Hom_R(R[r^{-1}],-) }   &&& \prod_{r \in R} R[r^{-1}] \Mod  \ar[lll]_{\operatorname{res}_r \circ \pi_r}      }.
$$
and the family $\{ (\Proj_r, R[r^{-1}] \Mod)  \}_{r \in R}$ of the trivial cotorsion pairs.
\begin{equation*}
\begin{split}
(\operatorname{res}_*( \Proj_*))^\perp
&=\{ N \in R\Mod \mid \mbox{Ext}_R^1(R[r^{-1}] ,N)=0 \textrm{ for every } r \in R   \}\\
&=\{ \textrm{contraadjusted $R$-modules} \}
\end{split}
\end{equation*}
$R$-modules in $^\perp((\operatorname{res}_*( \Proj_*))^\perp)$ are called very flat $R$-modules; see \cite[Section~1.1]{Pos17}

\section{Functor categories and Evaluation Functor}
In this and upcoming sections, we center on functor categories in order to obtain further results on cotorsion pairs induced by adjoint triples.

Throughout, $R$ denotes a ring with identity, and $\A$ denotes  a small preadditive category.   For any object $A$  in $\A$, we let $R_A$ denote  the endomorphism ring $\Hom(A,A)$. Note that   for any two  objects $A$ and $B$ in $\A$, $\Hom(A,B)$ has a canonical  $R_B$-$R_A$-bimodule structure induced from the composition rule
\begin{equation}\label{composition}
\xymatrix{
\Hom(B,B) \otimes_{\Z} \Hom(A,B) \otimes_{\Z} \Hom(A,A) \ar[rr] && \Hom(A,B)
}
\end{equation}

Since for every morphisms $g:A \rightarrow B$, $f:B \rightarrow B$ and $h:B\rightarrow C$ in $\mathcal{A}$, ${h \circ (f \circ g)=h \circ f \circ g=(h \circ f) \circ g}$,  the composition
is canonically factorized over $\Hom(B,C) \otimes_{R_B} \Hom(A,B)$
 $$\xymatrix{ \Hom(B,C)\otimes_{\Z} \Hom(A,B) \ar[rr] \ar[rd] && \Hom(A,C) \\
 & \Hom(B,C) \otimes_{R_B} \Hom(A,B) \ar[ru] &
 }
$$

\subsection{Evaluation Functor.}\label{eval_func} We let $\Add(\A,R\Mod)$ denote the category of additive functors from $\A$ to $R \Mod$. The evaluation functor
 $$\ev: \Add(\A,R\Mod) \otimes \A \rightarrow R \Mod$$ is the functor defined by $\ev(X,A)=X(A)=:\ev_A(X)$.

For a given additive functor $X:\A \longrightarrow R\Mod$ and objects $A$ and $B$ in $\A$, $X(A)$ has a canonical left $R_A$-module structure defined by $f \cdot a:=X(f)(a)$, where $f \in R_A$ and $a \in X(A)$. So it is easy to verify that  the associated homomorphism of abelian groups
$$X_{AB}:\ \xymatrix{ \Hom(A,B) \ar[r] & \Hom_R(X(A),X(B))}$$
induces the following left $R_B$ and $R_A$-module homomorphisms, respectively,
$$\overline{X}_{AB}: \Hom(A,B) \otimes_{ R_A} X(A) \longrightarrow X(B),  \quad g\otimes a \rightsquigarrow X(f)(x),$$
$$\underline{X}_{AB}: X(A) \longrightarrow \Hom_{R_B}(\Hom(A,B),X(B)),  \quad \underline{X}_{AB}( a)(g)= X(g)(a).$$

Besides,  as for every $f \in R_A$, $X(f)$ is an $R$-module homomorphism, we have
$$f \cdot (r \cdot a)=r \cdot (f \cdot a)$$ for every  $a \in X(A)$ and $r \in R$.

If $\alpha=\{\alpha_A\}_{A \in \A}:X \Rightarrow Y$ is a natural transformation of additive functors ${X,Y: \A \rightarrow R \Mod}$, then for any $f \in R_A$ and $a \in X(A)$, we have
$$
\alpha_A(f \cdot a)=\alpha_A(X(f)(a))=Y(f) ( \alpha_A(a))=f \cdot \alpha_A(a).
$$
In other words, for every object $A$ in $\A$,  $\alpha_A: X(A) \rightarrow Y(A)$ is an  $R_A$-module homomorphism. So
 the evaluation functor $\ev_A: \Add(\A,R \Mod) \rightarrow R \Mod$ has  also image in $R_A \Mod$
 \begin{equation}\label{eval_factor}
\xymatrix{ \Add( \A, R\Mod)\otimes \A \ar[rr]^-{\ev_A} \ar[rd]_{\ev_A} && R \Mod \\
 & R_A \Mod &
 }
  \end{equation}

\subsection{Remark.} Similarly,
 the evaluation functor $\ev_A: \Add(\A^ {\textrm{op}},R \Mod) \rightarrow R \Mod$ has   image in $R_A^{\textrm{op}} \Mod=\mbox{Mod-}R_A$.

 From now on, we assume that $\A$ satisfies the following condition.

\subsection{Hypothesis.}\label{hypot_A} For every object $A $ in  $\A$, there exists a ring morphism $\delta_A:R \rightarrow R_A$  such that for every morphism $f: A\rightarrow B$  in $\A$, the following diagram commutes
$$
\xymatrix{
&R_B\ar[rd]^{f^*}& \\
R  \ar[ru]^{\delta_B} \ar[rd]_{\delta_A}& & \Hom(A,B)\\
& R_A \ar[ru]_{f_*} &
}
$$

\vspace{4mm}
Hypothesis \ref{hypot_A} implies that for every object $A$ in $\A$,
$$\Hom(A,-) \in \Add(\A, R \Mod) \quad \textrm{  and } \quad  \Hom(-,A) \in \Add (\A^{op}, \mbox{Mod-}R).$$
Indeed, the $R_B$-$R_A$-bimodule $\Hom(A,B)$ is clearly $R$-$R$-bimodule by restriction of scalars. Besides, it satisfies the following equality
\begin{equation}\label{comm1}
r \cdot f= \delta_B(r) \circ f=f \circ \delta_A(r)=f \cdot r,
\end{equation}
for every  $r\in R$ and $f \in \Hom(A,B)$.
 Besides, for  any  morphisms $f:A \rightarrow B$ and  $g :B\rightarrow B'$ in $\A$ and an element $r$ in $R$, we have
$$g_*(r \cdot f)=g_*(f \circ \delta_A(r))=g \circ f \circ \delta_A(r)=\delta_{B'}(r) \circ (g \circ f)=r \cdot ( g_*(f)).$$
In other words, $g_*: \Hom(A,B) \rightarrow \Hom(A,B')$ is a left $R$-module homomorphism. Similarly, $g^*: \Hom(B',A) \rightarrow \Hom(B,A)$ is a right $R$-module homomorphism.

Furthermore, the evaluation functor $\ev_A$ given in \eqref{eval_factor} is factorized over the forgetful functor $U: R_A \Mod \rightarrow R \Mod$, and the composition rule
is canonically factorized as
 $$\xymatrix{ \Hom(B,C)\otimes_{R} \Hom(A,B) \ar[rr] \ar[rd] && \Hom(A,C)\\
 & \Hom(B,C) \otimes_{R_B} \Hom(A,B) \ar[ru]& \\
}$$

\subsection{Remark.} The conditions imposed on $\A$ in Hypothesis \eqref{hypot_A} is equivalent to being $\A$ an $R$-$R$-bimodule enriched category.

\subsection{Example.}\label{ex:chain1} We let $\A$ be the category whose set of objects is $\Z$, and
$$\Hom(i,j):=
\begin{cases}
R, & i=j \textrm{ or } i=j+1\\
0, & \textrm{ otherwise}
\end{cases}$$
It is easy to verify that $\Add(\A,R\Mod)$ is equivalent to the category   $ \bfC(R)$ of chain complexes of left $R$-modules.

\subsection{Example.}\label{ex:morita_context}
Consider a Morita context ring $\begin{bmatrix}
T & M\\
N & S
\end{bmatrix}$, that is,   $T$, $S$ are rings, $M$ and $N$ are $T$-$S$ and $S$-$T$-bimodules, respectively, together with  bimodule homomorphisms
$M \otimes_S N \rightarrow T$ , and $N \otimes_T M \rightarrow S$.

 We let $\A$ be the category with two objects $A$ and $B$ and  with morphisms $R_A:=T$, $R_B:=S$,  $\Hom(A,B):=N$, $\Hom(B,A):=M$ together with the canonical composition rule induced from the ring multiplications and bimodule homomorphisms of the given Morita context ring.
 It is immediate that
$$\Add(\A, \Ab) \cong \begin{bmatrix}
T & M\\
N & S
\end{bmatrix} \Mod. $$

\subsection{}
For any object $A$ in $\A$ and a  left $R_A$-module $M$, the authors in  \cite[Section~2]{CT13}  define the following $\Ab$-valued functors
\begin{align*}
{ \rm q}_A(M)(-) & :=\Hom(A,-) \otimes_{R_A} M\\
{ \rm p}_A(M)(-) &:=\Hom_{R_A}(\Hom(-,A),M),
\end{align*}
By the hypothesis  imposed on $\A$ in \eqref{hypot_A}, they are indeed $R \Mod$-valued functors, and hence, induce the functors
$${ \rm q}_A,{ \rm p}_A:  R_A \Mod \longrightarrow \Add(\A,R\Mod).$$


\subsection{Proposition.} For every object $A$ in $\A$,  $( { \rm q}_A, \ev_A, { \rm p}_A)$ is an adjoint triple.
\begin{proof}
By \cite[Proposition 1]{CT13}, viewing any functor $X \in \Add(\A, R \Mod)$ as an $\Ab$-valued functor, we already know that the natural transformation
$$\operatorname{Nat}({ \rm q}_A(M),  X ) \rightarrow \Hom_{R_A}(M,X(A)), \quad \alpha \rightsquigarrow  M \cong R_A \otimes_{R_A} M \xrightarrow{\alpha_A} X(A)$$
is injective. Surjectivity follows from the fact that the induced morphism
$$\overline{X}_{AB}: \Hom(A,B) \otimes_{R_A} X(A) \rightarrow  X(B)$$  is a left $R$-module homomorphism.

By dual arguments, it is easy to verify that $(\ev_A, { \rm p}_A)$ is an adjoint pair, as well.

\end{proof}

We let ${ \rm q}, { \rm p}:\prod_{A\in \A }R_A \Mod \rightarrow  \Add(\A,R\Mod) $ denote the functors $ \oplus_{A \in A} { \rm q}_A$ and $\prod_{A \in \A}{ \rm p}_A $, respectively.

\subsection{Proposition.}\label{prop:adj_ev2}
 $(  { \rm q}, ( \ev_A)_{A \in \A}, { \rm p} )$ is an adjoint triple.
 $$\xymatrix{ \prod_{A\in \A }R_A \Mod \ar@/^2em/[rrr]^{  { \rm q}} \ar@/_2em/[rrr]_{{ \rm p}} &&&    \Add(\A,R\Mod) \ar[lll]_{( \ev_A)_{A \in \A}} }.
$$

\subsection{Hypothesis.}\label{degreewise_exact}
From now on,  $ \prod_{A\in \A }R_A \Mod$ is assumed to be equipped with a fixed exact structure $ \prod_ {A \in \A} \mathcal{E}_A $, unless otherwise is stated. Accordingly, we assume that $\Add(\A, R\Mod)$ is equipped with the exact structure $\mathcal{E}$, which consists of  the   short exact sequences  $\mathbb{E}$  such that $\ev_A (\mathbb{E})$ belongs to $\mathcal{E}_A$ for every ${A \in \A}$. Clearly, $\mathcal{E}$ is a weakly idempotent complete exact structure on $\Add(\A, R\Mod)$. 

\subsection{Proposition.} \label{prop:generating}
For any object $A \in \mathcal{A}$, the following statements hold.
\begin{enumerate}[i)]
\item  $\ev_A$ is an exact functor for every object $A$ in $\A$. 
\item If $\prod_{A \in \A} \mcF_A$ is a generating class in
$\prod_{A \in \A}R_A \Mod$, then  so is $ \{ { \rm q}_A(\mcF_A) \}_{A \in \A}$ in $\Add(\A, R\Mod)$.
\item  If $\prod_{A \in \A} \mcG_A$ is a cogenerating class in
$\prod_{A \in \A}R_A \Mod$, then  so is $ \{ { \rm p}_A(\mcF_A) \}_{A \in \A}$ in $\Add(\A, R\Mod)$.

\end{enumerate}
\begin{proof}
\sloppy	The first statement  is immediate. It also implies that  for any  $\mathcal{E}_A$-projective object  $M$ in $R_A \Mod$,  ${ \rm q}_A(M)$ is projective in $\Add(\A, R\Mod)
$, as well.
	
	Now we prove the second statement. The dual arguments can be applied for the third statement. Let ${X \in \Add(\A, R\Mod)}$. For any object $A$  in $\A$, we consider an admissible epic ${F_A \twoheadrightarrow X(A)}$ in $R_A \Mod$  with $F_A \in \mcF_A$. By adjunction given in Proposition \ref{prop:adj_ev2},  there exists a natural transformation $\alpha: \oplus_{A \in \A} { \rm  q}_A(F_A) \Rightarrow X$. It is in fact an admissible epic in $\Add(\A, R \Mod)$. Indeed, for every ${B \in \A}$, the morphism $\alpha_B$  is the coproduct of the morphism
	$\oplus_{B \neq A \in \A} { \rm q}_A(F_A)(B) \longrightarrow X(B)$ and the admissible epic
	${ \rm q}_{B}(F_{B})(B) \cong F_{B} \twoheadrightarrow X(B)$ in $R_{B} \Mod$. By Proposition \ref{prop:sum_admissible}, $\alpha_B$ is an admissible epic in $R_{B} \Mod; \mathcal{E}_{B}$.
\end{proof}

As a particular case of Proposition \ref{prop:generating}, for any functor $X \in \Add(\A, R\Mod) $ the counit morphism $ { \rm q}((X(A))_{A \in \mathcal{A}}) \Rightarrow X$ and the unit morphism $ X \Rightarrow { \rm p}((X(A))_{A \in \mathcal{A}}) $ are an admissible epic and monic, respectively.

Applying  Corollary \ref{setup:3} for the adjoint triple given in Proposition \ref{prop:adj_ev2}, we obtain a more general form of Example \ref{ex:chain_ev_cot}. 
\subsection{Proposition.}\label{prop:cotorsion1}
Let $(\ \prod_{A \in \A} \mcF_A, \prod_{A \in \A} \mathcal{G}_A)$ be a  cotorsion pair in $ \prod_{A \in \A} R_A \Mod$.
\begin{enumerate}[i)]
\item
If for every $A \in \A$, any object in $\mcF_A$ is ${ \rm q}_A$-flat, then we have
$${ \rm q}(  \prod_{A \in \A} \mcF_A)^\perp=\{ X \in \Add(\mathcal{A}, R\Mod) \mid\ X(A) \in \mcG_A \textrm{ for every }A\in \A\}.$$
	\item
	If for every $A \in \A$, every object in $\mcG_A$ is ${ \rm p}_A$-coflat, then
	$$^\perp{ \rm p}(  \prod_{A \in \A} \mcG_A)=\{ X \in \Add(\mathcal{A}, R\Mod) \mid\ X(A) \in \mcF_A \textrm{ for every }A\in \A \}$$

	\end{enumerate}

\subsection{Example.} Consider the functor category representation of $\bfC(R)$ as given in Example \ref{ex:chain1}. For a given $i \in \mathbb{Z}$, the functors ${\rm q}_i$ and ${\rm p}_i$ are precisely the disc chain complexes $\operatorname{D}^i$ and $\operatorname{D}^{i+1}$ functors (see Example \ref{ex:chain_complex1}). Since $\operatorname{D}^i$ is an exact functor, the cotorsion pairs in Example \ref{ex:chain_ev_cot} are easily obtained  by  Proposition \ref{prop:cotorsion1}.

\vspace{3mm}

As a particular case, if we consider the absolute exact structure on $\prod_{A \in \A}  R_A \Mod$, then we have the following result.
\subsection{Corollary.}\label{cor:absol_eval}
Let $(\ \prod_{A \in \A} \mcF_A, \prod_{A \in \A} \mathcal{G}_A)$ be a  cotorsion pair in $\prod_{A \in \A} R_A \Mod$ equipped with the absolute exact structure.
\begin{enumerate}[i)]
\item
If  ${ \rm L_1 p}_A(\mcF_A)=\mbox{Tor}_{R_A}^1(\Hom(A,-),\mcF_A)=0$ for every object $A$ in $\A$, then we have
$${ \rm q}(  \prod_{A \in \A} \mcF_A)^\perp=\{ X \in \Add(\mathcal{A}, R\Mod) \mid\ X(A) \in \mcG_A  \textrm{ for every }A\in \A \}.$$

\item
If   ${ \rm R^1q}_A( \mcG_A)=\Ext^1_{R_A}( \Hom(-,A),\mcG_A)=0$ for every object  $A$ in $\A$, then
$$^\perp{ \rm p}(  \prod_{A \in \A} \mcG_A)=\{ X \in \Add(\mathcal{A}, R\Mod) \mid\ X(A) \in \mcF_A \textrm{ for every }A\in \A \}.$$	
	\end{enumerate}

 \subsection{Remark.} The categories $\prod_{A \in \A} R_A \Mod$ and $\Add(\mathcal{A}, R\Mod)$ equipped with the absolute exact structures have enough injectives and projectives. Under the conditions given in   Corollary \ref{cor:absol_eval}-(i), if $(\ \prod_{A \in \A} \mcF_A, \prod_{A \in \A} \mathcal{G}_A)$ is hereditary, then clearly the class ${ \rm q}(  \prod_{A \in \A} \mcF_A)^\perp$ is closed under quotients of monomorphism, and therefore the cotorsion pair generated by ${ \rm q}(  \prod_{A \in \A} \mcF_A)$ is hereditary, as well. By the same reason, under the conditions given in   Corollary \ref{cor:absol_eval}-(ii), the cotorsion pair cogenerated by $\perp{ \rm p}(  \prod_{A \in \A} \mcG_A)$ is hereditary, as well.

\section{Stalk functor in functor categories}
 We continue studying on $\Add(\A, R \Mod)$. Unless otherwise stated, we still assume that $\Add(\A, R \Mod)$ is equipped with the exact structure $\mathcal{E}$ induced from a fixed exact structure $ \prod_{A \in \A} \mathcal{E}_A$ on $\prod_{A \in \A }R_A \Mod $ as stated in \eqref{degreewise_exact}.  In this section,  we investigate cotorsion pairs arising through `stalk functor'. In order  to do it, we need to impose further conditions on $\A$.

\subsection{Hypothesis.}\label{hyp2}  We assume that $\A$ satisfies the following condition: For every two objects $A\neq B$ in $\A$, the composition rule
$$\Hom(A,B) \otimes_{R_A} \Hom(B,A) \longrightarrow \Hom(B,B)$$
is zero.

\subsection{Example.} Consider a Morita context ring as given in  Example \ref{ex:morita_context}. Whenever the bimodule homomorphisms
 $M \otimes_S N \rightarrow T$ and $N \otimes_T M \rightarrow S$ are zero, it is  called \textit{Morita context ring with zero trace ideals}.  In particular, when $M=0$ or $N=0$, then it is called \textit{a formal
triangular matrix ring}.

The associated category $\A$ to a Morita context ring with zero trace ideals as presented in Example \ref{ex:morita_context} satisfies Hypothesis \ref{hyp2}.

\subsection{Stalk functor.}\label{def:stalk} For  every objects  $A$ and $B$ in $\A$ and $M \in R_A \Mod$, we define
$$
{ \rm s}_A(M)(B):=
\begin{cases}
0; & A \neq B\\
M; & A=B
\end{cases}
$$
For a given morphism $f: B \rightarrow B'$ in $\A$, we define a morphism
$${ \rm s}_A(M)(f): { \rm s}_A(M)(B) \rightarrow { \rm s}_A(M)(B')$$
 as follows:  If $B \neq A$ or $B' \neq A$, then ${ \rm s}_A(M)(f):=0$; if $B=A=B'$,  then ${ \rm s}_A(M)(f):=f\cdot(-)$, the left $R_A$-module scalar multiplication on $M$.

\subsection{Proposition.}\label{prop:stalk1}
For every object $A$ in $\A$, ${ \rm s}_A: R_A \Mod \longrightarrow \Add(\A, R\Mod)$ is a functor.

\begin{proof}
The only thing that we need to show is that for every left $R_A$-module $M$, ${ \rm s}_A(M)(-)$ belongs to  $\Add(\A, R \Mod)$.

It is immediate that for every $B \in \A$, ${ \rm s}_A(M)(B)$ is a left $R$-module.  By \eqref{comm1}, for a given morphism $f\in R
_A $, the morphism ${ \rm s}_A(M)(f)$ is an $R$-module homomorphism.

Clearly, for every $B\in B$, ${ \rm s}_A(M)(\id_B)=\id$.  It is left to prove that ${ \rm s}_A(M)(-)$ preserves compositions. Consider two arbitrary morphisms $f:B \rightarrow B'$ and $g: B' \rightarrow B''$ in $\A$. It is immediate that
$${ \rm s}_A(M)(g \circ f) ={ \rm s}_A(M)(g) \circ { \rm s}_A(M)(f) $$ whenever  $B \neq A$ or $B'' \neq A$. Suppose that $B=A=B''$. If $B'=A$,  from the associativity of the scalar multiplication,
$${ \rm s}_A(M)(g \circ f):=(g \circ f)\cdot (-)=(g \cdot f)(-)=g \cdot (f \cdot-)={ \rm s}_A(M)(g) \circ { \rm s}_A(M)(f).$$
If $B' \neq A$, by Hypothesis \ref{hyp2}, $g \circ f=0$, and hence,
 $${ \rm s}_A(M)(g \circ f)=0= { \rm s}_A(M)(g) \circ { \rm s}_A(M)(f).$$
	
\end{proof}	

\subsection{Remark.}\label{rem:difference_stalk}
The stalk functor introduced in \eqref{def:stalk} differs subtly from the one defined in \cite[Definition 7.9]{HJ22}. For instance, let $\A$ be the additive category generated by the quiver   
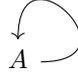
\vspace{-8mm}
$$
\begin{tikzcd}
A \arrow[out=0,in=90,loop]
\end{tikzcd}
$$
Then $\A$ has only one object with the endomorphism ring $R_A=\Z[x]$. For every polynomial $p(x) \in \Z[x]$, we have $p(x)=p(0)+\overline{p}(x) \in \Z \oplus \langle x \rangle=\Z[x]$. We let $S_A: \Ab \rightarrow \Add(\A, \Ab)$ denote the functor defined in \cite[Proposition 7.15]{HJ22}. 
$$
\xymatrix{\Z[x]\Mod \ar[r]^-{{\rm s}_A} & \Add(\Z, \Ab) & Ab \ar[l]_-{S_A}}
$$ 
If $ M \in \Z[x] \Mod$, then $S_A(M)(p(x))=p(0)\cdot (-)\neq { \rm s}_A(M)(p(x))$.
\vspace{5mm}

In \eqref{phi_psi},  we introduce several crucial constructions which will be used both for the  left/right of ${ \rm s}_A$,  and to induce certain short exact sequences  in Proposition \ref{Prop:left_Derived_coeq}.

\subsection{}\label{phi_psi}
Let $X: \A \longrightarrow R \Mod$ be  an additive functor, and $A$ be an object in $\A$.
By the universal property of coproducts, the family $\{ \overline{X}_{BA}\}_{ \substack{ B \in \A \\ A \neq B}}$ of $R_A$-module homomorphisms (see \eqref{eval_func}) induces a unique $R_A$-module homomorphism
$$\xymatrix{ \varphi_A^X: \displaystyle \bigoplus_{ \substack{ B \in \A \\ A \neq B}} \Hom(B,A) \otimes_{R_B} X(B) \ar[rr] && X(A)}.$$
We let ${ \rm c}_A(X):=\coker (\varphi_A^X)$. By the universal property of coproducts and cokernels, it in fact leads to   a functor
$$ \xymatrix{ { \rm c}_A: \Add(\A, R\Mod) \ar[rr] && R_A\Mod}. $$

On the other hand, there exist two canonical $R_A$-module homomorphisms $h^1_{A,X}$ and $h^2_{A,X}$ making the following  upper and lower square commutative, respectively,
$$
\scalebox{0.85}{
\xymatrix{
\Hom(B',A) \otimes_{R_{B'}} \Hom(B,B') \otimes_{R_B} X(B) \ar[rr]^-{\id \otimes \overline{X}_{BB'}} \ar@{^{(}->}[d] && \Hom(B',A) \otimes_{R_{B'}} X(B')  \ar@{_{(}->}[d] \\
\displaystyle \bigoplus_{ \substack{ B,B' \in \A \\ A \neq B,B'}} \Hom(B',A) \otimes_{R_{B'}} \Hom(B,B') \otimes_{R_B} X(B) \ar@<-1ex>[rr]_-{h^2_{A,X}} \ar@<+1ex>[rr]^-{h^1_{A,X}} && \displaystyle \bigoplus_{ \substack{ B \in \A \\ A \neq B}} \Hom(B,A) \otimes_{R_A} X(B)\\
	\Hom(B',A) \otimes_{R_{B'}} \Hom (B,B') \otimes_{R_B} X(B)\ar[rr]_-{ \circ \otimes \id } \ar@{_{(}->}[u] && \Hom(B,A) \otimes_{R_{B}} X(B)  \ar@{^{(}->}[u]}}$$
where $\circ$ denotes the composition rule. Note that the objects $B$ and $B'$ in the above diagrams can be assumed to be different.  Besides, we let
${\nu_{A,X}:=(  \id \otimes \overline{X}_{AB})_{ \substack{ B \in \A \\ A \neq B}}}$
$$
\scalebox{0.85}{
\xymatrix{
\Hom(B,A) \otimes_{R_{B}} \Hom(A,B) \otimes_{R_A} X(A) \ar[rr]^-{\id \otimes \overline{X}_{AB}} \ar@{^{(}->}[d] && \Hom(B,A) \otimes_{R_{B}} X(B)  \ar@{^{(}->}[d] \\
\displaystyle \bigoplus_{ \substack{ B \in \A \\ A \neq B}} \Hom(B,A) \otimes_{R_{B}} \Hom(A,B) \otimes_{R_A} X(A) \ar[rr]_-{\nu_{A,X}}   && \displaystyle \bigoplus_{ \substack{ B \in \A \\ A \neq B}} \Hom(B,A) \otimes_{R_A} X(B)}}$$

As $\A$ satisfies Hypothesis \ref{hyp2}, and     $X$ is a functor,
we have the following commutative diagrams
\begin{equation}\label{diag4}
\scalebox{0.9}{
\xymatrix{
	\Hom(B,A) \otimes_{R_{B}} \Hom (A,B) \otimes_{R_A} X(A) \ar[rr]^-{\id \otimes \overline{X}_{AB}} \ar[d]&& \Hom(B,A) \otimes_{R_{B}} X(B) \ar[d]^{\overline{X}_{BA}} \\
0 \ar[rr]&& X(A)
}}
\end{equation}

\begin{equation}\label{diag3}
\scalebox{0.9}{\xymatrix{
	\Hom(B',A) \otimes_{R_{B'}} \Hom (B,B') \otimes_{R_B} X(B) \ar[rr]^-{\id \otimes \overline{X}_{BB'}} \ar[d]_{\circ \otimes \id}&& \Hom(B',A) \otimes_{R_{B'}} X(B') \ar[d]^{\overline{X}_{B'A}} \\
	\Hom(B,A) \otimes_{R_B} X(B) \ar[rr]_{\overline{X}_{BA}}&& X(A)
}}
\end{equation}
for every objects $A \neq B,B'$ in $\A$, and therefore,  the morphism $\varphi_A^X$ is factorized over the cokernel of  the homomorphism $\phi_A^X:=(h^1_{A,X}-h^2_{A,X}, \nu_{A,X})$

$$\xymatrix{
\displaystyle\bigoplus_{ \substack{ B \in \A \\ A \neq B}} \Hom(B,A) \otimes_{R_B} X(B)  \ar[rr]^-{ \varphi_A^X} \ar[rd]_-{\sigma_A} && X(A)\\
 &\coker(\phi_A^X) \ar[ru] &
}$$
where $\sigma_A$ denotes the canonical quotient homomorphism. Note that
$$ \coker(\phi_A^X)  \cong ( \bigoplus_{A \neq B \in \A} \Hom(B,A) \otimes_{R_B} X(B) )/ <g \circ f \otimes x -g \otimes X(f)(x), s \otimes X(l)(y)>$$
where $x \in X(B)$, $y \in X(A)$ and  $f:B \rightarrow B'$, $g: B' \rightarrow A$, $s: B \rightarrow A$, $l: A \rightarrow B$ are any morphisms in $\A$ with $B\neq B'$ and $B,B' \neq A$.

Similarly, applying universal property of products for the family $\{ \underline{X}_{AB}\}_{ \substack{ B \in \A \\ B \neq A}}$ of $R_A$-module homomorphisms, we have the canonical left $R_A$-module homomorphism
$$\psi_A^X: X(A) \longrightarrow \prod_{A \neq B \in \A}\Hom_{R_B}(\Hom(A,B),X(B)).$$
 We let ${ \rm k}_A(X):= \ker (\psi_A^X)$, which   leads to a functor
$${ \rm k}_A: \Add(\A,R \Mod) \longrightarrow R_A \Mod.$$

Again,  there exist two canonical $R_A$-module homomorphisms $t^1_{A,X}$ and $t^2_{A,X}$ making the following  upper and lower square commutative, respectively,

$$
\scalebox{0.8}{
\xymatrix{
\Hom_{R_B}(\Hom(A,B),X(B)) \ar[rr]^-{ {\underline{X}_{BB'}}_*} && \Hom_{R_B}( \Hom(A,B), \Hom_{R_{B'}}( \Hom(B,B'), X(B')))  \\
\displaystyle\prod_{A \neq B \in \A}\Hom_{R_B}(\Hom(A,B),X(B)) \ar@<-1ex>[rr]_-{t^2_{A,X}} \ar@<+1ex>[rr]^-{t^1_{A,X}} \ar@{->>}[u]  \ar@{->>}[d]  &&
\displaystyle\prod_{A \neq B,B' \in \A}\Hom_{R_B}( \Hom(A,B), \Hom_{R_{B'}}( \Hom(B,B'), X(B'))) \ar@{->>}[u]  \ar@{->>}[d]  \\
	\Hom_{R_B}(\Hom(A,B),X(B)) \ar[rr]_-{\divideontimes}  && \Hom_{R_{B'}}( \Hom(A,B'), \Hom_{R_{B}}( \Hom(B',B), X(B)))  }}$$
where $\divideontimes$ denotes the composition of
$$ \xymatrix{ \Hom_{R_B}(\Hom(A,B),X(B)) \ar[r]^-{\circ^*} &  \Hom_{R_{B}}( \Hom(B',B) \otimes_{R_{B'}}\Hom(A,B'),  X(B))}$$
 with the canonical isomorphism
$$\scalebox{0.8}{ \xymatrix{ \Hom_{R_{B}}( \Hom(B',B) \otimes_{R_{B'}}\Hom(A,B'),  X(B)) \ar[r]^-{ \cong }&  \Hom_{R_{B'}}( \Hom(A,B'), \Hom_{R_{B}}( \Hom(B',B), X(B))) }}.$$
Besides, consider the $R_A$-module homomorphism $\tau_{A,X}:= \prod_{A \neq B \in \A} {\underline{X}_{BA}}_*$
$$
\scalebox{0.8}{
\xymatrix{
\displaystyle\prod_{A \neq B \in \A}\Hom_{R_B}(\Hom(A,B),X(B))  \ar[rr]^-{\tau_{A,X}}  \ar@{->>}[d]  &&
\displaystyle\prod_{A \neq B \in \A}\Hom_{R_B}( \Hom(A,B), \Hom_{R_{A}}( \Hom(B,A), X(A)))   \ar@{->>}[d]  \\
	\Hom_{R_B}(\Hom(A,B),X(B)) \ar[rr]_-{{\underline{X}_{BA}}_*}  && \Hom_{R_{B}}( \Hom(A,B), \Hom_{R_{A}}( \Hom(B,A), X(A)))  }}$$

It is easy to verify that the $R_A$-module homomorphism $\psi_A^X$ is factorized over the kernel  of $\mu_{A,X}:=(t^1_{A,X}-t^2_{A,X}, \tau_{A,X})$
$$\xymatrix{
  X(A)  \ar[rr]^-{ \psi_A^X} \ar[rd] && \displaystyle  \prod_{A \neq B \in \A}\Hom_{R_B}(\Hom(A,B),X(B))\\
 &\ker (\mu_{A,X}) \ar@{^{(}->}[ru] &
}$$
\subsection{Proposition.}\label{prop:adj_stalk} For any object  $A$ in $\A$, $( { \rm c}_A, { \rm s}_A, { \rm k}_A)$ is an adjoint triple
$$\xymatrix{ \Add(\A, R \Mod) \ar@/^2em/[rrr]^{{ \rm c}_A} \ar@/_2em/[rrr]_{{ \rm k}_A} &&&      R_A \Mod  \ar[lll]_{{ \rm s}_A} }.
$$

\begin{proof}
 We construct a natural  isomorphism
$$\gamma: \xymatrix{\Hom_{R_A}({ \rm c}_A(-), -) \ar@{=>}[r] & \mbox{Nat}_R(-, { \rm s}_A(-))}.$$
Let $z  \in \Hom_{R_A}( { \rm c}_A(X), M)  $. We let  $\gamma(z)_A: X(A) \longrightarrow M= { \rm s}_A(M)(A)$ be the composition of $z$ with the canonical quotient homomorphism $ X(A) \rightarrow { \rm c}_A(X)$, and $\gamma(z)_B:=0: X(B) \rightarrow 0={ \rm s}_A(M)(B)$ for $B \neq A$. Clearly, $\gamma(z)_B$ is an $R$-module homomorphism for every object $B$ in $\A$. Now we show that ${\gamma(z):=\{\gamma(z)_B\}_{B \in \A}}$ is a natural transformation from $X$ to ${ \rm s}_A(M)$. For it, for a given morphism ${f: B \rightarrow B'}$ in $\A$, we need to show the commutativity of  the following diagram
$$
\xymatrix{X(B) \ar[rr]^{X(f)} \ar[d]_{\gamma(z)_B}&& X(B') \ar[d]^{\gamma(z)_{B'}} \\
{ \rm s}_A(M)(B) \ar[rr]_{ { \rm s}_A(M)(f)} && { \rm s}_A(M)(B')	
}
$$
It is immediate whenever $B' \neq A$. If $B =A=B'$, then  it follows from the fact that $\gamma(z)_A$ is an $R_A$-module homomorphism. Now suppose that $B \neq A$ and $B' = A$. $X(f)$ is factorized as
$$\xymatrix{X(B) \ar[r]^-{f \otimes \id} &  \Hom(B,A) \otimes_{R_A} XB  \ar@{^{(}->}[r] & \bigoplus_{A \neq B \in \A} \Hom(B,A) \otimes_{R_B} XB \ar[r]& XA}.$$
Therefore, $\gamma(\alpha)_A \circ X(f)=0.$ As a consequence, $\gamma: X \Rightarrow { \rm s}_A(M)$ is a natural transformation.  One can easily show the naturality and injectivity  of $\gamma$.

Let $\alpha: X \Rightarrow { \rm s}_A(M)$ be a natural transformation. As pointed out  in \eqref{eval_func},  $\alpha_A$ is an $R_A$-module homomorphism. Besides, for every $B \neq A$, the following diagram commutes
$$
\xymatrix{\Hom(B,A) \otimes_{R_B} X(B) \ar[rr]^-{\overline{X}_{AB}} \ar[d]_{\id \otimes \alpha_B}&& X(A) \ar[d]^{\alpha_{A}} \\
	0=\Hom(B,A) \otimes_{R_B} s_A(M)(B) \ar[rr] && M	
}
$$
So $\alpha_A$ is factorized over $c_A(X)$, which shows the surjectivity of $\gamma$.

Hence, $( { \rm c}_A, { \rm s}_A)$ is an adjoint pair. Similar arguments can be applied for $( { \rm s}_A,{ \rm k}_A)$.

\end{proof}

 We let
 $${ \rm c}, { \rm k}:\Add(\A, R \Mod) \longrightarrow \prod_{A \in \A} R_A \Mod, $$
  $${{ \rm s}:\prod_{A \in \A} R_A \Mod \longrightarrow  \Add(\A, R \Mod)}$$ denote the functors $({ \rm c}_A)_{A \in \A},({ \rm k}_A)_{A \in \A} $ and   $\oplus_{A \in \A} { \rm s}_A $, respectively.
\subsection{Proposition.}\label{prop:adj2_stalk} $( { \rm c},  { \rm s}, { \rm k})$ is an adjoint triple
$$\xymatrix{ \Add(\A, R \Mod) \ar@/^2em/[rrr]^{{ \rm c}} \ar@/_2em/[rrr]_{{ \rm k}} &&&     \prod_{A \in \A} R_A \Mod  \ar[lll]_{ { \rm s}} }.
$$
\begin{proof}
Note that 	for every object $B$ in $\A$,
$$\bigoplus_{A \in \A} { \rm s}_A(M_A)(B) = { \rm s}_B(M_B)= \prod_{A \in \A}{ \rm s}_A(M_A)(B).$$
Therefore,
$${ \rm s}((M_A)_{A \in \A})=\bigoplus_{A \in \A} { \rm s}_A(M_A)=\prod_{A \in \A}{ \rm s}_A(M_A).$$
So we have
\begin{equation*}
\begin{split}
\Hom_{\prod_{A \in \A}  R_A \Mod}(({ \rm c}_A(X))_{A \in \A}, (M_A)_{A \in \A}) & := \prod_{A \in \A} \Hom_{R_A}( { \rm c}_A(X), M_A)\\
& \cong \prod_{A \in \A } \mbox{Nat}_R(X, { \rm s}_A(M_A))\\
& \cong \mbox{Nat}_R(X, \prod_{A \in \A} { \rm s}_A(M_A))\\
\end{split}
\end{equation*}
Similarly,  $ ({ \rm s}, { \rm k})$ is proved to be an adjoint pair.
\end{proof}

\subsection{Remark.} \sloppy It is immediate that  the functor ${ \rm s}$ is an exact functor from $\prod_{A \in \A} R_A \Mod$ to $ \Add( \A,  R\Mod) $.
However,   neither ${ \rm c}_A$ nor ${ \rm k}_A $  preserve short exact sequences in $\Ext^1(-, { \rm s}_A(M))$ and $\Ext^1({ \rm s}_A(M),-)$, respectively, as  required in Corollary \ref{setup:3}-(iii-iv).  The category of  chain complexes of left $R$-modules, which is a very particular case of $\Add(\A, R\Mod)$, can be considered as a counterexample (see Example  \ref{ex:chain_stalk2}).  For a more general case, for example,  assume that there exist two different  objects $A$ and $B$ in $\A$ such that $\Hom(A,B)$ is a flat right $R_A$-module. For any non-trivial left $R_A$-module $M$, consider the short exact sequence
$$\xymatrix{ K_{M,A} \ar@{^{(}->}[r] & { \rm q}_A(M) \ar@{->>}[r] & { \rm s}_A(M) }$$
 given in Proposition \ref{prop:ses1}. Note that
the associated $R_A$-module homomorphism
$$\underline{ { \rm q}_A(M)}_{AB}:
\xymatrix{ { \rm q}_A(M)(A)\cong M \ar[r] & \Hom_{R_B}( \Hom(A,B),{ \rm q}_A(M)(B) )},$$
${ \rm q}_A(M)(A)(m)=- \otimes m: \Hom(A,B) \rightarrow \Hom(A,B)\otimes_{R_A} M$, is injective, and therefore, ${ \rm k}_A( { \rm q}_A(M))=0$. Applying ${ \rm k}_A$ to the aforementioned short exact exact sequence,  we  have $0 \rightarrow 0 \rightarrow M$.

\vspace{3mm}

Applying  Corollary \ref{setup:3} to the adjoint triples given in Proposition \ref{prop:adj_stalk} and Proposition \ref{prop:adj2_stalk}, we obtain the following two corollaries. 

\subsection{Corollary.}\label{corollary:stal-cogen-cot.} Let  $A$ be an object in  $\A$. For a given   cotorsion pair $(\mcF_A,  \mcG_A)$ in $R_A \Mod$, the following statements  hold:
\begin{enumerate}[(i)]
\item If  $\mcF_A$ is a generating class in $R_A \Mod$, then
$${ \rm s}_A(  \mcF_A)^\perp= \biggl \{ X: \A \rightarrow R \Mod \biggm| \
\begin{array}{l}
{ \rm k}_A(X) \in \mcG_A \textrm{  and  }\\
$X$\textrm{ is a }{ \rm k}_A \textrm{-coflat object}
\end{array}
 \biggr\} .$$

\item If  $\mcG_A$ is a cogenerating class in $R_A \Mod$, then
$$^\perp{ \rm s}_A(  \mcG_A)=\biggl \{ X: \A \rightarrow R \Mod \biggm| \
\begin{array}{l}
{ \rm c}_A(X) \in \mcF_A \textrm{  and  }\\
$X$\textrm{ is a }{ \rm c}_A \textrm{-flat object}
\end{array}
 \biggr\}
.$$
\end{enumerate}

\subsection{Corollary.} For any   cotorsion pair $(  \prod_{A \in \A} \mcF_A,  \prod_{A \in \A} \mcG_A)$ in $\prod_{A \in \A} R_A \Mod$, the following statements  hold:
\begin{enumerate}[(i)]
\item If  $\prod_{A \in \A} \mcF_A$ is a generating class in $\prod_{A \in \A} R_A \Mod$, then
$${ \rm s}( \prod_{A \in \A} \mcF_A)^\perp=\biggl \{ X: \A \rightarrow R \Mod \biggm| \
\begin{array}{l}
{ \rm k}_A(X) \in \mcG_A \textrm{  and  }$X$\textrm{ is a }{ \rm k}_A \textrm{-coflat}\\
 \textrm{object for every } A \in \A
\end{array}
 \biggr\}.$$

\item If  $\prod_{A \in \A} \mcG_A$ is a cogenerating class in $\prod_{A \in \A} R_A \Mod$, then
$$^\perp{ \rm s}( \prod_{A \in \A}  \mcG_A)=\biggl \{ X: \A \rightarrow R \Mod \biggm| \
\begin{array}{l}
{ \rm c}_A(X) \in \mcF_A \textrm{  and  }$X$\textrm{ is a }{ \rm c}_A \textrm{-flat}\\
 \textrm{object for every } A \in \A
\end{array}
 \biggr\}.$$
\end{enumerate}

\subsection{Remark.} From Lemma \ref{lemma:ext_adjoint_pair}-(i-ii), it is immediate that ${\rm p}_B(G) \in { \rm s}_A(  \mcF_A)^\perp$, and ${\rm q}_B(F) \in {}^\perp{ \rm s}_A(  \mcG_A)$ for every $A,B \in \A$, $F\in \mcF_B$ and $G \in \mathcal{G}_B$. By Proposition~\ref{prop:generating}, if $\prod_{A \in \A} \mcG_A$ is a cogenerating class in $\prod_{A \in \A} R_A \Mod$, so is ${ \rm s}( \prod_{A \in \A} \mcF_A)^\perp$; if $\prod_{A \in \A} \mcF_A$ is a generating class in $\prod_{A \in \A} R_A \Mod$, so is $^\perp{ \rm s}(  \prod_{A \in \A} \mcG_A)$.

 \subsection{Remark.} If we consider the absolute exact structure on both categories
$\Add(\A, R \Mod)$ and $\prod_{A \in \A} R_A \Mod$, then by  Proposition~\ref{prop:adjoints} and Proposition~\ref{lemma:derived_zero}, we have
\begin{align*}
{ \rm s}( \prod_{A \in \A} \mcF_A)^\perp&= \{ X: \A \rightarrow R \Mod \mid\ { \rm k}_A(X) \in \mcG_A \textrm{  and $ { \rm R^1k}_A(X)=0$   for every $A \in \A$} \},\\
^\perp{ \rm s}( \prod_{A \in \A}  \mcG_A)&=\{ X: \A \rightarrow R \Mod \mid\ { \rm c}_A(X) \in \mcF_A \textrm{  and $ { \rm L_1c}_A(X)=0$   for every $A \in \A$} \}\\
\end{align*}
which are equal to the classes $ \Psi(\prod_{A \in \A}  \mcG_A )$    and $ \Phi(\prod_{A \in \A}  \mcF_A)$, respectively, introduced in \cite[Definition 1.5,Theorem 1.7]{HJ19}.

\vspace{5mm}

The notion of (co)flatness related to a functor is, in general, abstract and not easily visualized. To better understand ${\rm c}_A$-flat and ${\rm k}_A$-coflat objects in $\Add(\A, R\Mod)$, we aim to describe them in terms of specific short exact sequences, similar to the one given in Example \ref{ex:chain_stalk2}. To do this, we first need to prove the following auxiliary results.

\subsection{Lemma.}\label{prop:ses1} For every object $A$ in $\A$ and  a left $R_A$-module $M$,  the canonical natural transformations 
$${ \rm q}_A(M) \Longrightarrow  { \rm s}_A(M)\quad \textrm{ and } \quad   { \rm s}_A(M)  \Longrightarrow { \rm p}_A(M)$$
are admissible epic and admissible monic, respectively, in $\Add(\A, R\Mod)$. 

\begin{proof}
We only prove that  the first morphism is an admissible epic as the argument for the  second one is similar. Let $\alpha: { \rm q}_A(M) \Rightarrow { \rm s}_A(M)$ be the counit of the  adjunction $({ \rm q_A}, { \rm ev_A})$ for the functor ${ \rm s}_A(M)$. It is easy to verify that
$$\alpha_A:R_A \otimes_{R_A} M \rightarrow M$$ is the canonical isomorphism induced from the  left scalar multiplication of $M$ by elements in $R_A$, and $\alpha_B=0: \Hom(A,B) \otimes_{R_A} M \rightarrow 0$  for $B \neq A$.  So for every object $B$ in $\A$, the morphism $\alpha_B$ is an admissible epic in $\mathcal{E}_B$. By construction of $\mathcal{E}$ (see  \eqref{degreewise_exact}), $\alpha$ is an admissible epic in $\Add(\A, R\Mod)$. 
 \end{proof}

\subsection{Remark.} The canonical morphisms  given in Lemma \ref{prop:ses1} induce  the following canonical short exact sequences in $\Add(\A,R\Mod)$
$$\xymatrix{ { \rm K}_{M,A} \ar@{^{(}->}[r] & { \rm q}_A(M) \ar@{->>}[r] & { \rm s}_A(M) }$$
$$ \xymatrix{ { \rm s}_A(M) \ar@{^{(}->}[r] & { \rm p}_A(M) \ar@{->>}[r] & { \rm C}_{M,A}   }  $$
From the construction, it is immediate that the functors ${ \rm K}_{M,A}, { \rm C}_{M,A}: \A \longrightarrow R \Mod$ satisfies
$${ \rm K}_{M,A}(B):=\begin{cases}
0 &; B=A\\
{ \rm q}_A(M)(B)& ;B \neq A
\end{cases}
$$
$${ \rm C}_{M,A}(B):=\begin{cases}
0 & ;B=A\\
{ \rm p}_A(M)(B) &; B \neq A
\end{cases}
$$

\subsection{Lemma.}\label{lemma:coker_C_M}
For any  functor $X \in \Add( \A, R\Mod)$ and a left $R_A$-module $M$, we have
\begin{align*}
\mbox{Nat}_R(X, { \rm C}_{M,A}) &\cong \Hom_{R_A}(\coker (\phi_{A,X}) ,M);\\
\mbox{Nat}_R( { \rm K}_{M,A}, X) & \cong \Hom_{R_A}(M, \ker (\mu_{A,X}) ).
\end{align*}

\begin{proof}
Due to the duality of arguments, we only prove the isomorphism
$\mbox{Nat}_R(X, { \rm C}_{M,A}) \cong \Hom_{R_A}(\coker (\phi_{A,X}) ,M)$.  Let $\alpha:X \Rightarrow { \rm C}_{M,A} $ be a natural transformation. From \eqref{eval_func}, we already know that $\alpha_B$ is an $R_B$-module homomorphism for any object $B$. Now, for a given object $B$ in $\A$ with $B \neq A$,  we let $\overline{\alpha}_B$ denote the $R_A$-module homomorphisms in $\Hom_{R_A}( \Hom(B,A) \otimes_{R_B} X(B),M  ) \cong \Hom_{R_B}(X(B), { \rm p}_A(M)(B))$ corresponding to $\alpha_B$ under the usual adjunction. Note that $\overline{\alpha}_B(f \otimes x)=\alpha_B(x)(f) $. It is easy to verify that  for every two objects $B,B'$ in $\A$ different from $A$,  we have
$$\overline{\alpha}_{B'} \circ ( \id \otimes \overline{X}_{BB'})= \overline{\alpha}_B \circ ( \circ \otimes \id)$$
as given in the diagram \eqref{diag3}, and  since $\alpha$ is a natural transformation,  ${\alpha_B \circ X(l)=0}$ for every morphism $l: A \rightarrow B$, and hence, $\overline{\alpha}_B \circ (\id \otimes \overline{X}_{AB})=0$ (see the diagram~\eqref{diag3}). So there exists a unique $R_A$-module homomorphism $z_\alpha: \coker (\phi_{A,X}) \longrightarrow M$ satisfying
$$z_\alpha \circ \sigma_A= \Sigma_{A \neq B} \overline{\alpha}_B: \displaystyle\bigoplus_{A \neq B \in \A} \Hom(B,A) \otimes_{R_B} X(B) \longrightarrow M. $$
Clearly, the assignment
$$\mbox{Nat}_R(X, { \rm C}_{M,A}) \longrightarrow \Hom_{R_A}(\coker (\phi_{A,X}) ,M), \quad \alpha \rightsquigarrow z_\alpha, $$  is well defined, and is natural in both $X$ and $M$. Since $\sigma_A$ is an epimorphism, it is also injective.

Now, let $z:\coker (\phi_{A,X}) \longrightarrow M $ be an $R_A$-module homomorphism. We let $\alpha_A:=0: X(A) \longrightarrow 0$, and for an object $B \neq A$, we let $\alpha_B$ denote the $R_B$-module  homomorphism in $ \Hom_{R_B}(X(B), { \rm p}_A(M)(B)) \cong \Hom_{R_A}( \Hom(B,A) \otimes_{R_B} X(B),M  )    $  corresponding to the $R_A$-module homomorphism $z \circ \sigma_A $ restricted to $\Hom(B,A) \otimes_{R_B} X(B)$ under the usual adjunction. Note that $\alpha_B(x)=z \circ \sigma_A(- \otimes x).$ By restriction of scalars, the family ${\alpha:=\{\alpha_B \}_{B \in \A}}$ is clearly  a family of $R$-module homomorphisms. We claim that $\alpha$ is a natural transformation from $X$ to $ { \rm C}_{M,A}$. For it,  consider a morphism $f: B \rightarrow B'$ in $\A$.  The case when  $B'=A$ is immediate, as  $\alpha_A \circ X(f)=0= { \rm C}_{M,A}(f) \circ \alpha_B$. Now suppose that $B' \neq A$. Then there are two cases to consider: $B \neq A$ and $B =A$.

\underline{$B \neq A$}: For any morphism $g:B' \longrightarrow A$ and $x \in X(B)$, we have
\begin{align*}
({ \rm C}_{M,A}(f) \circ \alpha_B) (x)(g)&= (z \circ \sigma_A)((g \circ f) \otimes x)= (z \circ \sigma_A)(g \otimes X(f)(x))\\
&=\alpha_{B'}(X(f)(x))(g)=(\alpha_{B'} \circ  X(f))(x)(g),
\end{align*}
which implies ${ \rm C}_{M,A}(f) \circ \alpha_B=\alpha_{B'} \circ  X(f) $.

\underline{$B=A$}: For any object $x$ in $X(A)$, we have
$$(\alpha_{B'} \circ X(f))(x)(-)=z \circ \sigma_A(-\otimes X(f)(x)): \Hom(B',A) \longrightarrow M$$
But $\sigma_A(g\otimes X(f)(x))=0$ for any morphism $g \in \Hom(B',A)$. So
$$\alpha_{B'} \circ X(f)=0 = { \rm C}_{M,A}(f) \circ \alpha_A.$$

As a consequence, $\alpha$ is a natural transformation, and therefore, the assignment $\mbox{Nat}_R(X, { \rm C}_{M,A}) \longrightarrow \Hom_{R_A}(\phi_{A,X}) ,M)$ defined above is surjective.
\end{proof}

\subsection{Lemma.}\label{lemma:compos_q_c} For every two objects $A \neq B$ in $\A$,
\begin{enumerate}[(i)]
\item ${\rm c}_A \circ { \rm q}_A=\id$ and  ${ \rm c}_A \circ { \rm q}_B=0$;
\item ${\rm k}_A \circ { \rm p}_A=\id$ and   ${ \rm k}_A \circ { \rm p}_B=0$.
\end{enumerate}
\begin{proof}
By adjunctions,  we have
\begin{align*}
\Hom_{R_A}({ \rm c}_A ( { \rm q}_A(-)),- )& \cong \mbox{Nat}_R({ \rm q}_A(-), { \rm s}_A(-) )\\
& \cong \Hom_{R_A}(-,{ \rm ev}_A ( { \rm s}_A(-) ))\\
&= \Hom_{R_A}(-,\id(-) ).
\end{align*}
$$\Hom_{R_A}({ \rm c}_A ( { \rm q}_B(-)),- ) \cong \mbox{Nat}_R({ \rm q}_B(-), { \rm s}_A(-) ) \cong \Hom_{R_B}(-,{ \rm ev}_B ( { \rm s}_A(-) ))=0,$$
So ${ \rm c}_A  \circ  { \rm q}_A=\id$ and ${ \rm c}_A \circ { \rm q}_B=0$. By similar arguments, we have the equalities given in (ii).
\end{proof}

Now we can provide a concrete description of ${\rm c}_A$-flat and ${\rm k}_A$-coflat functors, which will also be crucial in proving Theorem \ref{them:proj-inj}. Recall morphisms $\phi_{A,X}$ and $\ker (\mu_{A,X}) $ introduced in  \eqref{phi_psi}.
\subsection{Proposition.}\label{Prop:left_Derived_coeq} Let $A$ be an object in $\A$. For any additive functor  ${X: \A \longrightarrow R \Mod} $, the following statements hold.
\begin{enumerate}[(i)]
\item $X$   is ${ \rm c}_A$-flat in $\Add(\A, R\Mod)$ if and only if the canonical sequence
	$$\coker (\phi_{A,X})\longrightarrow X(A) \longrightarrow { \rm c}_A(X)$$
	given in \eqref{phi_psi} is a short exact sequence  in $ R_A \Mod$.
	
\item
$X$   is ${ \rm k}_A$-coflat in $\Add(\A, R\Mod)$ if and only if the canonical sequence
	$${ \rm k}_A(X) \longrightarrow  X(A) \longrightarrow    \ker (\mu_{A,X})  $$
	given in \eqref{phi_psi} is a short exact sequence  in $ R_A \Mod$.
	
\end{enumerate}	
\begin{proof}
We only prove the first statement.

$(\Rightarrow)$ Suppose that $X$ is a ${ \rm c}_A$-flat object in $\Add(\A, R\Mod)$.  We firstly show that the morphism
$ X(A) \rightarrow { \rm c}_A(X)$ is an admissible epic in $R_A \Mod$. For it, using the identity morphism $\id_{X(B)}$ for every object $B$ in $\A$ and the adjoint pair $( { \rm q}, ({\rm ev}_B)_{B \in \A})$, we have  the canonical natural transformation $$\alpha:\bigoplus_{B\in \A} { \rm q}_B(X(B)) \Rightarrow X,$$
which is an admissible epic in $\Add(\A, R\Mod)$ just as pointed out in the proof of Proposition \ref{prop:generating}. By assumption, ${\rm c}_A(\alpha)$ is an admissible epic in $R_A\Mod$. As ${ \rm c}_A$ is a left adjoint functor,  it preserves coproducts, and using Lemma \ref{lemma:compos_q_c}-(i), we have
$${\rm c}_A(  \bigoplus_{B\in \A} { \rm q}_B(X(B))) \cong \bigoplus_{B\in \A} {\rm c}_A(   { \rm q}_B(X(B)) \cong X(A). $$
So ${\rm c}_A(\alpha)$ is isomorphic to  the canonical quotient morphism $X(A) \rightarrow {\rm c}_A(X)$. In other words, it is an admissible epic in $R_A \Mod$.

On the other hand, let  $E$ be an injective cogenerator left $R_A$-module in the absolute sense  ( that is, an injective cogenerating object in $R_A \Mod$ regarding  to the absolute exact structure). It is clearly an  injective object  with respect to the exact structure $\mathcal{E}_A$, as well.   Consider the following canonical short exact sequence in $\Add(\A, R\Mod)$ as given in Lemma \ref{prop:ses1}
$$ \mathbb{E}: \xymatrix{ { \rm s}_A(E) \ar@{^{(}->}[r] & { \rm p}_A(E) \ar@{->>}[r] & { \rm C}_{E,A}   }. $$
Applying the functor $\mbox{Nat}_R(X,-)$, it induces the following long exact sequence of abelian groups
$$\scalebox{0.75}{\xymatrix{ 0 \ar[r]& \mbox{Nat}_R(X,{ \rm s}_A(E)) \ar[r] & \mbox{Nat}_R(X,{ \rm p}_A(E)) \ar[r] & \mbox{Nat}_R(X,{ \rm C}_{E,A}) \ar[r]& \Ext^1_{\mathcal{E}}(X,{ \rm s}_A(E)) \ar[r] & \Ext^1_{\mathcal{E}}(X,{ \rm p}_A(E))  }}.$$
Since ${ \rm ev}_A$ is an exact functor, $ { \rm p}_A(E)$ is an injective object in $\Add(\A, R\Mod)$. So $\Ext^1_{\mathcal{E}}(X,{ \rm p}_A(E))=0$.
Besides, as  $X$ is a ${\rm c}_A$-flat object in $\Add(\A, R\Mod)$, by Proposition \ref{prop:adjoints}, we have
 $\Ext_{\mathcal{E}}^1(X,  { \rm s}_A(E))\cong \Ext^1_{\mathcal{E}_A}({ \rm c}_A( X), E)=0$, and therefore, $\mbox{Nat}_R(X, \mathbb{E})$ is a short exact sequence of abelian groups. By adjunctions and Lemma~\ref{lemma:coker_C_M}, it is isomorphic to
$$\scalebox{0.9}{ \xymatrix{ 0 \ar[r] & \Hom_{R_A}( { \rm c}_A(X), E) \ar[r] & \Hom_{R_A}(X(A), E) \ar[r] &  \Hom_{R_A}(  \coker(\phi_{A,X}), E) \ar[r] & 0}.} $$
As $E$ is a cogenerating left $R_A$-module, $ \coker(\phi_{A,X})$ is the kernel of the admissible epic $X(A) \rightarrow {\rm c}_A(X)$, and therefore,
$$\xymatrix{ \coker (\phi_{A,X}) \ar@{^{(}->}[r] &  X(A) \ar@{->>}[r] & { \rm c}_A(X)}$$
is a short exact sequence in $R_A\Mod$.

$(\Leftarrow)$ Suppose that the canonical sequence
$$\xymatrix{ \coker (\phi_{A,X}) \ar@{^{(}->}[r] &  X(A) \ar@{->>}[r] & { \rm c}_A(X)}$$
is a short exact sequence in $R_A \Mod$. By Lemma \ref{lemma:lem}-(i), it is enough to show that ${ \rm c}_A$ preserves short exact sequences in $\Ext^1(X, { \rm s}_A(- ))$. Let $\mathbb{E}$ be a short exact sequence in $\Add(\A, R\Mod)$ of the form
$$\xymatrix{ { \rm s}_A(M ) \ar@{^{(}->}[r] &  Z \ar@{->>}[r]^-\alpha & X}, $$
where $M$ is a left $R_A$-module.
Without lost of generality, we can assume that $X(B)=Z(B)$ for every object $B$ in $\A$ with $A \neq B$. So we have $h^1_{A,X}=h^1_{A,Z}$ and $h^2_{A,X}=h^2_{A,Z}$. On the other hand, we have the following commutative diagram
$$
\scalebox{0.8}{\xymatrix{
\displaystyle \bigoplus_{A\neq B \in \A} \Hom(B,A)\otimes_{R_B} \Hom(A,B) \otimes_{R_A} Z(A) \ar[rr]^{\oplus \id \otimes \id \otimes \alpha_A} \ar[d]_{\nu_{A,Z}} && \displaystyle \bigoplus_{A\neq B \in \A} \Hom(B,A)\otimes_{R_B} \Hom(A,B) \otimes_{R_A} X(A)\ar[d]^{\nu_{A,X}}\\
\displaystyle \bigoplus_{A\neq B \in \A} \Hom(B,A)\otimes_{R_B}  Z(B) \ar@{=}[rr] && \displaystyle \bigoplus_{A\neq B \in \A} \Hom(B,A)\otimes_{R_B}  X(B)
}}
$$
whose upper horizontal arrow  is an epimorphism. So the images of $\nu_{A,X}$ and $\nu_{A,Z}$ coincide, and  therefore, we have $\coker(\phi_{ A,X})=\coker(\phi_{A,Z})$
together with the following commutative diagram
$$
\scalebox{0.75}{\xymatrix{
0 \ar@{^{(}->}[r] \ar[dd]  & \displaystyle \bigoplus_{A\neq B \in \A} \Hom(B,A)\otimes_{R_B}Z(B) \ar@{=}[rrr]\ar[dd]_{\varphi_A^Z} \ar[rd] &&&  \displaystyle  \bigoplus_{A\neq B \in \A} \Hom(B,A)\otimes_{R_B}X(B)  \ar[dd]^{\varphi_A^X}  \ar[ld] \\
&&\coker(\phi_{ A,Z}) \ar[ld] \ar@{=}[r]&\coker(\phi_{ A,X}) \ar@{^{(}->}[rd] &\\
M \ar@{^{(}->}[r]  &Z(A)  \ar@{->>}[rrr]&&& X(A)\\
}}
$$
By \cite[Proposition 2.16]{Buh}, the canonical morphism $\coker(\phi_{ A,Z}) \longrightarrow Z(A)$ is an admissible monic in $R_A$, as well. Applying Snake Lemma (see \cite[Corollary 8.13]{Buh}) to the following commutative diagram

$$
\xymatrix{
0 \ar[d] \ar[r]&\coker(\phi_{ A,Z}) \ar@{^{(}->}[d] \ar@{=}[r]&\coker(\phi_{ A,X}) \ar@{^{(}->}[d]\\
M \ar@{^{(}->}[r]  &Z(A)  \ar@{->>}[r]& X(A)\\
}
$$
we obtain the following short exact sequence in $R_A \Mod$
$$\xymatrix{ M={\rm c}_A( { \rm s}_A(M)) \ar@{^{(}->}[r] &  {\rm c}_A(Z) \ar@{->>}[r] & { \rm c}_A(X)}.$$

\end{proof}

\subsection{Theorem.}\label{cor:cot_final} For any   cotorsion pair $(  \prod_{A \in \A} \mcF_A,  \prod_{A \in \A} \mcG_A)$ in $\prod_{A \in \A} R_A \Mod$, the following statements  hold.
\begin{enumerate}[(i)]
\item If  $\prod_{A \in \A} \mcF_A$ is a generating class in $\prod_{A \in \A} R_A \Mod$, then
$${ \rm s}( \prod_{A \in \A} \mcF_A)^\perp=\biggl\{ X: \A \rightarrow R \Mod \biggm \vert \
\begin{array}{l}
{ \rm k}_A(X) \in \mcG_A,
X(A) \longrightarrow    \ker (\mu_{A,X})\\
\textrm{is an admissible epic in } R_A\Mod,\ \forall A \in \A
\end{array} \biggr \}.$$

\item If  $\prod_{A \in \A} \mcG_A$ is a cogenerating class in $\prod_{A \in \A} R_A \Mod$, then
$$^\perp{ \rm s}( \prod_{A \in \A}  \mcG_A)=\biggl \{ X: \A \rightarrow R \Mod \biggm \vert \
\begin{array}{l}
 { \rm c}_A(X) \in \mcF_A,  \coker (\phi_{A,X}) \hookrightarrow X(A)\\
\textrm{is an admissible monic in  } R_A\Mod,\ \forall A \in \A
\end{array}
\biggr \}.$$
\end{enumerate}

\section{Projective and injective functors}
In this section, we apply Theorem~\ref{cor:cot_final} to provide an intrinsic characterization of projective and injective functors. We  retain all the hypotheses assumed on  $\mathcal{A}$ so far. Applying Theorem~\ref{cor:cot_final}-(i) to  the trivial cotorsion pair $(   \prod_{A \in \A} R_A \Mod, \prod_{A \in \A} \mathcal{E}_A\mbox{-Inj} )$, we have  induces 
$${ \rm s}( \prod_{A \in \A} R_A \Mod)^\perp=\left \{ X: \A \rightarrow R \Mod \biggm \vert \
\begin{array}{l}
{ \rm k}_A(X) \in \mathcal{E}_A \mbox{-Inj},\\
X(A) \longrightarrow    \ker (\mu_{A,X}) \textrm{ is an admissible}\\
\textrm{epic in } R_A\Mod,\  \forall A \in \A
\end{array} \right \}.$$

Similarly, by Theorem~\ref{cor:cot_final}-(ii), the trivial cotorsion pair $( \prod_{A \in \A} \mathcal{E}_A\mbox{-Proj},  \prod_{A \in \A} R_A \Mod)$ leads to 
$$^\perp{ \rm s}( \prod_{A \in \A}  R_A \Mod)=\left \{ X: \A \rightarrow R \Mod \biggm \vert \
\begin{array}{l}
 { \rm c}_A(X) \in \mathcal{E}_A \mbox{-Proj},\\
   \coker (\phi_{A,X}) \hookrightarrow X(A) \textrm{ is an admissible}\\
\textrm{monic in  } R_A\Mod,\  \forall A \in \A
\end{array}
\right \}.$$

The following containments in $\Add(\A, R\Mod)$ can be easily verified
$$ \mathcal{E} \mbox{-Inj} \subseteq  { \rm s}( \prod_{A \in \A}  R_A \Mod)^\perp \quad \textrm{ and } \quad    \mathcal{E} \mbox{-Proj} \subseteq  {}^\perp{ \rm s}( \prod_{A \in \A}  R_A \Mod). $$
The  classes on the right side of the above inclusions   have a precise description and  are easier  to work from a homological point of view. Furthermore, in the very well-known cases of functor categories (e.g. the category of chain complexes of left $R$-modules, the module category over a formal triangular matrix ring etc.), they give rise to injective and projective objects in functor categories. We need to impose further conditions on $\mathcal{A}$ to obtain the reverse containments.

\subsection{Condition.}\label{condition1} There doesn't exist a chain of morphisms in $\A$
$$\xymatrix{\cdots \ar[r] & A_3 \ar[r]^{f_2} & A_2 \ar[r]^{f_1} & A_1}$$
satisfying $A_i \neq A_j$ for every $1 \leq i \neq j$, and $f_1 \circ \cdots \circ f_n \neq 0$ for every $n \geq 1$.

\subsection{Condition.}\label{condition2} There doesn't exist a chain of morphisms in $\A$
$$\xymatrix{A_1 \ar[r]^{f_1}  & A_2 \ar[r]^{f_2} & A_3 \ar[r]^{f_3} & \cdots }$$
satisfying $A_i \neq A_j$ for every $1 \leq i \neq j$, and $f_n \circ \cdots \circ f_1 \neq 0$ for every $n \geq 1$. 


\subsection{Lemma.}\label{lemma:coker_zero}
For any functor $X \in \Add( \A, R \Mod)$, the following statements hold.
\begin{enumerate}[(i)]
\item If Condition \ref{condition1} is satisfied, then  $X=0$ if and only if ${ \rm c}_A(X)=0$ for every $A \in \A$.
\item If  Condition \ref{condition2} is satisfied, then $X=0$ if and only if ${ \rm k}_A(X)=0$ for every $A \in \A$.
\end{enumerate}
\begin{proof}
The necessity part is clear. Now suppose that ${ \rm c}_A(X)=0$ for every $A \in \A$. It implies that the $R_A$-module homomorphism  $\varphi_A^X$ is an epimorphism for every object $A$ in $\A$. For the sake of contradiction, we assume that $X \neq 0$. Consider an object  $A_1$  in $\A$ such that $X(A_1) \neq 0$. Let $ x_1$ be a non-zero object  in $X(A_1)$.
By assumption,  there exist a finite set $I_2$ of objects in $\A$ different from $A_1$, and a family $\{(f_B,x_B)\}_{B \in I_2}$, where $f_B \in \Hom(B,A_1)$ and $x_B \in X(B)$, such that
$$\sum_{B\in I_2}X(f_B)(x_B)=x_1.$$
Without lost of generality, we can assume that $X(f_B)(x_B)\neq 0$ for every $B \in I_2$. We choose any object object in $I_2$, and we denote it by $A_2$, and $f_1:=f_{A_2}$, $x_2:=x_{A_2}$. Similarly, as $x_2 \neq 0$, there exist a finite family $I_3$ of objects in $\A$ different from $A_2$, and a family $\{(f_B,x_B)\}_{B\in I_3}$, where $f_B \in \Hom(B,A_2)$ and $x_B \in X(B)$, such that
$$\sum_{B\in I_3}X(f_B)(x_B)=x_2.$$
We claim that there is an object $B_0\in I_3$ such that
$$X(f_1 \circ f_{B_0})(x_{B_0})=(X(f_1) \circ X(f_{B_0}))(x_{B_0})\neq 0.$$
Indeed, if $X(f_1 \circ f_{B})(x_{B})=0$ for every $B \in I_3$, then   we would have
$$X(f_1)(x_2)=X(f_1)(  \sum_{B\in I_3} X(f_{B})(x_{B}) )=\sum_{B\in I_3} X(f_1 \circ f_B)(x_B)=0,$$
which is a contradiction. We let $A_3$ be  any object in $I_3$ such that
$$(X(f_1 \circ f_{A_3}))(x_{A_3})\neq 0.$$
 We let $f_2:=f_{A_3}$ and $x_3:=x_{A_3}$. Note that $A_3 \neq A_2$. From  Hypothesis \ref{hyp2}, we also have  $A_3 \neq A_1$.

 By induction, we can construct a sequence of morphisms in $\A$
$$ \xymatrix{\cdots \ar[r] & A_3 \ar[r]^{f_2} & A_2 \ar[r]^{f_1} & A_1}  $$
with $A_i \neq A_j$  for every $i \neq j \geq1$ such that $X(f_1 \circ \ldots \circ f_n) \neq 0$, and therefore,
 $f_1 \circ \ldots \circ f_n \neq 0$ for every $n \geq 1$, which is a contradiction by assumption.  So $X=0$.
\end{proof}

The previous lemma is slightly different and more general than the result given in \cite[Theorem 7.19]{HJ22} as we have the following example
\subsection{Example.}\label{ex:non-left rooted}Let $\A$ be a category with objects $\mathbb{N}$ and morphisms
$$\Hom(n,m):=\begin{cases}
R & n \in \{m,\ldots, 2m\}\\
0 & \textrm{otherwise}
\end{cases}$$
It clearly satisfies  Condition \ref{condition1}, but it is not left rooted nor fulfills the conditions given  in \cite[Theorem 7.29]{HJ22}.

\vspace{4mm}

Following a   similar strategy used in the proof of \cite[Theorem 7.29]{HJ22}, we have the follwing characterization of projective and injective additive functors.
\subsection{Theorem.}\label{them:proj-inj}
\begin{enumerate}[(i)]
\item If Condition \ref{condition1} is satisfied, then  
$$\mathcal{E} \mbox{-Proj} = {}^\perp{ \rm s}( \prod_{A \in \A}  R_A \Mod).$$

\item If Condition \ref{condition2} is satisfied, then  
$$\mathcal{E} \mbox{-Inj} = { \rm s}( \prod_{A \in \A}  R_A \Mod)^\perp$$
\end{enumerate}
\begin{proof}
We only prove the first statement. 	We already know the inclusion $\subseteq$. Now, let $X \in{}^\perp{ \rm s}( \prod_{A \in \A}  R_A \Mod)$. Then for every $A \in \A$, we have the following canonical short exact sequence
\begin{equation}	\label{s.e.s}
	\xymatrix{ \coker(\phi_{A,X}) \ar@{^{(}->}[r] & XA \ar@{->>}[r]^-{\pi_A} & { \rm c}_A(X)    }
\end{equation}
 with an $\mathcal{E}_A$-projective  left $R_A$-module ${ \rm c}_A(X)$. So for every $A \in \A$,  there exists an admissible monic $\iota_A: { \rm c}_A(X) \hookrightarrow X(A)$ such that $\pi_A \circ \iota_A=\id$. 
 Using the counit of the adjunction $({ \rm q}, ({ \rm ev}_A)_{A \in \A})$, we obtain a natural transformation  $\alpha: \oplus_{A \in \A} { \rm q}_A( { \rm c}_A(X)) \Rightarrow X$. Note that $\alpha$ is the composition of  the splitting admissible monic $\oplus_{A \in \A} { \rm q}_A( \iota_A):  \oplus_{A \in \A} { \rm q}_A( { \rm c}_A(X)) \Rightarrow  \oplus_{A \in \A} { \rm q}_A( X(A))   $ with the canonical natural transformation $$\oplus_{A \in \A} { \rm q}_A(  X(A)) \Rightarrow X,  $$
 which is an admissible epic in $\Add(\A,R\Mod)$ as pointed out in the proof of Proposition \ref{prop:generating}. Note that $ \oplus_{A \in \A} { \rm q}_A( { \rm c}_A(X))$ is an $\mathcal{E}$-projective object in $\Add(\A, R\Mod)$

 By Lemma \ref{lemma:compos_q_c}, for any object  $B$ in $\A$, we have
 $${ \rm c}_B(\oplus_{A \in \A} { \rm q}_A( { \rm c}_A(X))  ) \cong \oplus_{A \in \A} { \rm c}_B( { \rm q}_A( { \rm c}_A(X)) \cong { \rm c}_A(X),$$
and therefore,  ${ \rm c}_B(X)$ is an isomorphism. Since
 ${ \rm c}_B$ preserves all colimits, ${ \rm c}_B( \coker(\alpha))=\coker({ \rm c}_B(\alpha) )=0$ for every object $B $ in $\A$. By Lemma~\ref{lemma:coker_zero}-(i), $\coker(\alpha)=0$, that is,  $\alpha$ is an epimorphism. 
 
Notice that as $R_A \Mod$ is an abelian category,  the short exact sequence \ref{s.e.s} in $(R_A \Mod; \mathcal{E}_A)$ is  a usual short exact sequence in $R_A \Mod$ with respect to the absolute exact structure. So by Proposition \ref{Prop:left_Derived_coeq}, for every object $A$ in $\A$, $X$ is a ${ \rm c}_A$-flat object whenever both categories $\Add(\A, R\Mod)$ and $R_A \Mod$ are considered to be equipped with the absolute exact structures. It implies that
$${ \rm c}_A( \ker(\alpha))=\ker({ \rm c}_A(\alpha) )=0 $$
for every object $A$ in $\A$, and hence, by  Lemma~\ref{lemma:coker_zero}-(i), $\ker(\alpha)=0$. As a consecuence, $\alpha$ is an isomorphism.
\end{proof}

\subsection{Remark.} If one starts with an exact structure $\mathcal{E}'$ on $R \Mod$, then by restriction of scalars, it leads to a family $\{ \mathcal{E}_A \}_{A \in \A}$ of exact structures on $R_A \Mod$ defined as follows: a short exact sequence $\mathbb{E}_A$ in $R_A \Mod$ belongs to $\mathcal{E}_A$ if and only if it belongs to $\mathcal{E}'$ as a short exact sequence of $R$-modules. Note that $\mathcal{E}_A$-projective and injective objects are of the form $R_A \otimes_R M$ and $\Hom_R (R_A, N)$, respectively, where $M \in \mathcal{E}'\mbox{-Proj}$ and $N \in \mathcal{E}'\mbox{-Inj}$.


\section*{Funding}
The first and the third named authors  were partially supported by grant I+D+i PID2020-113206GB-I00, funded by  MICIU/AEI/10.1303/501100011033 and by grant 22004/PI/22 funded by Fundaci\'on S\'eneca-Agencia de Ciencia y Tecnolog\'ia de la Regi\'on de Murcia.

The second named author was supported by grant I+D+i PID2020-113206GB-
I00, funded by MCIU/AEI/10.13039/501100011033.

\end{document}